\documentclass[review,onefignum,onetabnum]{siamart190516}
\usepackage{tikz,subfigure}
\usepackage{booktabs}
\usepackage[symbol]{footmisc}
\usepackage{algorithmic}
\usepackage{amsmath,amssymb}
\usepackage{esvect}
\usepackage{arcs}
\usepackage{multirow}
\usepackage{accents}
\usepackage{graphicx}

\newcommand{\Tau}{\mathrm{T}}
\newcommand{\z}{\hskip -.03in}
\newcommand{\REMOVE}[1]{}
\DeclareMathOperator*{\argmin}{arg\,min}
\DeclareMathAccent{\wtilde}{\mathord}{largesymbols}{"65}
\DeclareMathAccent{\wtildea}{\mathord}{largesymbols}{"65}
\makeatletter
\def\widebreve{\mathpalette\wide@breve}
\def\wide@breve#1#2{\sbox\z@{$#1#2$}%
     \mathop{\vbox{\m@th\ialign{##\crcr
\kern0.08em\brevefill#1{0.8\wd\z@}\crcr\noalign{\nointerlineskip}%
                    $\hss#1#2\hss$\crcr}}}\limits}
\def\brevefill#1#2{$\m@th\sbox\tw@{$#1($}%
  \hss\resizebox{#2}{\wd\tw@}{\rotatebox[origin=c]{90}{\upshape(}}\hss$}
\makeatletter
\def\linefill{%
\leavevmode
\leaders\hrule\hskip\dimexpr\textwidth-2pt\mbox{}}

\definecolor{color1}{RGB}{255,174,120}

\usepackage{lipsum}
\usepackage{amsfonts}
\usepackage{graphicx}
\usepackage{epstopdf}
\usepackage{algorithmic}
\ifpdf
  \DeclareGraphicsExtensions{.eps,.pdf,.png,.jpg}
\else
  \DeclareGraphicsExtensions{.eps}
\fi

\crefname{hypothesis}{Hypothesis}{Hypotheses}

\headers{Energy minimization H-curl AMG}{Tuminaro/Glusa}

\title{A structure preserving {H}-curl algebraic multigrid method for the eddy current equations\thanks{%Submitted to the editors DATE.
The authors were supported by the U.S.~Department of Energy, Office of Science, Office of Advanced Scientific Computing Research, Applied Mathematics program.  Sandia National Laboratories is a multimission laboratory managed and operated by National Technology and Engineering Solutions of Sandia, LLC., a wholly owned subsidiary of Honeywell International, Inc., for the U.S. Department of Energy's National Nuclear Security Administration under grant~DE-NA-0003525.  This paper describes objective technical results and analysis.  Any subjective views or opinions that might be expressed in the paper do not necessarily represent the views of the U.S. Department of Energy or the United States Government.}}%}

\author{
  Raymond S. Tuminaro\thanks{Computational Mathematics Department, Sandia National Laboratories, Livermore, CA
  (\email{rstumin@sandia.gov}, \url{http://www.cs.sandia.gov/ccr-rstumin}).}
\and
  Christian Glusa\thanks{Scalable Algorithms Department, Sandia National Laboratories, Albuquerque, NM
  (\email{caglusa@sandia.gov}). %, \url{http://www.cs.sandia.gov/ccr-rstumin}).
  }
}

\usepackage{amsopn}

\ifpdf
\hypersetup{
  pdftitle={A structure preserving {H}curl algebraic multigrid method for the eddy current equations } 
  pdfauthor={R. Tuminaro/C. Glusa}
}
\fi

\begin{document}

\maketitle

\begin{abstract}
A new algebraic multigrid method (AMG) is presented for solving the linear systems associated with the eddy current
approximation to Maxwell’s equations.  This AMG method extends an idea proposed by 
Reitzinger and Sch\"oberl. The main feature of the Reitzinger and Sch\"oberl algorithm (RSAMG) is that it maintains
null-space properties of the $\nabla \times \nabla \times$ operator throughout all levels of the AMG hierarchy.
It does this by enforcing a commuting relationship involving grid transfers and the discrete gradient operator.
This null-space preservation property is critical to the algorithm's success, however enforcing 
this commuting relationship is non-trivial except in the special case where one leverages a piece-wise constant nodal 
interpolation operator.  For this reason, mesh independent convergence rates are generally 
not observed for RSAMG due to its reliance on  sub-optimal piece-wise constant interpolation.
We present a new AMG algorithm that enforces the same commuting relationship. The main advance
is that the new structure preserving H-curl algorithm (SpHcurlAMG) does not rely on piece-wise constant interpolation and 
can leverage fairly general and more sophisticated nodal interpolation operators.
The key idea %behind SpHcurlAMG 
is to employ energy minimization AMG (EAMG) to construct edge interpolation grid
transfers and to enforce the commuting relationship by embedding it as constraints within an EAMG procedure.
While it might appear that solving such a constrained energy minimization is costly, we illustrate how this is 
not the case in our context.  Numerical results are then given demonstrating mesh independent convergence over a range of
test problems.
\end{abstract}



%%%%%%%%%%%%%%%%%%%%%%%%%%%%%%%%%%%%%%%%%%%%%%%%%%%%%%%%%%%%%%%%%%%%%
%%%%%%%%%%%%%%%%%%%%%%%%%%%%%%%%%%%%%%%%%%%%%%%%%%%%%%%%%%%%%%%%%%%%%
%%%%%%%%%%%%%%%%%%%%%%%%%%%%%%%%%%%%%%%%%%%%%%%%%%%%%%%%%%%%%%%%%%%%%
\section{Introduction}\label{sec:intro}
We consider the solution of linear systems associated with the 3D eddy current formulation of Maxwell's equations
given by
\begin{equation} \label{eq: eddy current}
\nabla \times \nabla \times \vec{u} + \sigma(x,y,z) ~ \vec{u} = f
\end{equation}
where
$\vec{u}$ is the electric field, $\sigma$ is the electrical conductivity, and $f$ is
the right hand side.  Many simulations require linear solutions to discrete versions of this partial differential equation (PDE) when modeling systems that include Maxwell's equations or
magnetohydrodynamic formulations, often as sub-calculations within the larger simulation.
Discretization by first-order edge elements on arbitrary
unstructured meshes gives rise to the discrete linear system
\begin{equation} \label{eq: discrete eddy}
A^{(e)}_h u_h = f_h
\end{equation}
that must be solved for $u_h$ where
%Throughout this paper,
the superscript $(e)$ denotes operators that are
applied to edge unknowns.
% while a superscript $(n)$ denotes an operator that is applied to node unknown (e.g., those associated with linear nodal finite elements).
The matrix $ A^{(e)}_h $ can be split into two pieces such that
$ A^{(e)}_h  = S^{(e)}_h + M^{(e)}_h $ where $S^{(e)}_h$ is the discretization of the
first term in \eqref{eq: eddy current} and $M^{(e)}_h $ is the discrete representation
of the second term.

It is well-known that a first-order edge element discretization preserves the
following vector identity
$$
\nabla \times ( \nabla \phi ) = 0
$$
which implies that the following
discrete relationship holds
\begin{equation} \label{eq: discrete null space condition}
S^{(e)}_h D^{(n\rightarrow e)}_h = {\Theta}
\end{equation}
where
$\Theta$ is the zero matrix.
$D^{(n\rightarrow e)}_h$ is a rectangular matrix that represents the discrete gradient operator.
The $n\rightarrow e$ superscript indicates a transformation between nodal and edge quantities.
$D^{(n\rightarrow e)}_h $ is easy to construct as each row
corresponds to a mesh edge and contains at most two nonzeros $( \pm 1 )$ each associated with an edge end point.
% in the %graph
%associated with the
%underlying mesh.
Rows with only one nonzero correspond to edges
where one node lies on a boundary associated with a Dirichlet condition. Structure preserving discretizations
that preserve the curl-gradient relationship have many advantages, but the large near null space
of $A^{(e)}_e$ implies that the problem is extremely ill-conditioned when $\sigma$ is small. For this
reason, general purpose iterative methods typically fail, and so special purpose algorithms must be considered.

Several multigrid approaches have been proposed for the eddy current PDE as well as
for related PDEs with similar large near null space properties. These include both
geometric and algebraic multigrid algorithms~\cite{
Arnold.D;Falk.R;winther.R2000,
https://doi.org/10.1002/nla.577,
BoGaHuRoTu03,
BoHuRoTu2003,
BoHuSiTu2008,
Hiptmair.R1997,
Hiptmair98,
ReSc02,
HiXu06,
HuTuBoGaRo06,
JoLe2005,
KoPaVa08,
doi:10.1137/110859361,
d27a9d32-2d8a-31ac-b344-f536a25e37d7}.
Most (though not all) successful algebraic multigrid (AMG) solvers for H-curl fall into either of two categories: structure preserving methods or auxiliary-based solvers.
The approach described in this manuscript is a structure preserving method and is most closely related to~\cite{ReSc02, BoGaHuRoTu03, HuTuBoGaRo06}. These methods (as well as most H-curl geometric multigrid methods)
 are structure preserving in the sense that the null space property \eqref{eq: discrete null space condition} is satisfied by the discretization matrix on each level of the multigrid hierarchy.
Generally, the AMG structure preserving methods are centered around AMG ideas that generate specialized grid transfer operators to guarantee this structure preserving property.
The original structure preserving AMG method was proposed in ~\cite{ReSc02} but this method is not scalable due to its reliance on on a sub-optimal piecewise constant nodal interpolation operator.
The algorithms in ~\cite{BoGaHuRoTu03, HuTuBoGaRo06} are extensions of ~\cite{ReSc02} that leverage smoothed aggregation ideas to allow for more sophisticated grid transfer operators.
However, these methods cannot be easily adapted to cases where weak connections must be internally dropped to address large material variations or anisotropic problems. When dropping occurs, the prolongator smoothing step must be modified to avoid a large increase in the number of matrix nonzeros on coarse grids while still assuring that the range space of interpolation includes the null space.  Unfortunately, this is not generally possible for the 
smoothed aggregation extension of H-curl.

The auxiliary H-curl AMG methods~\cite{HiXu06,KoPaVa08,doi:10.1137/110859361,d27a9d32-2d8a-31ac-b344-f536a25e37d7,BoHuSiTu2008} take a different approach that centers on developing 
preconditioners for a related problem that includes a vector Laplacian sub-system and a scalar Laplacian sub-system. The key advantage is that this preconditioner is also
effective for the desired H-curl problem and that it can leverage standard algebraic multigrid solvers for the sub-systems.  In some sense the auxiliary method allows one to apply an AMG algorithm that by itself would not be 
appropriate for H-curl problems. However, the standard AMG method can be applied to a reformulated problem so that it is effective for the H-curl problem when used
in conjunction with an additional AMG subproblem (which essentially corrects for deficiencies).  Though asymptotically optimal, the underlying problem reformulation does incur 
some costs (e.g., setup, storage and multiple V cycle invocations in the cycle ) and introduces some algorithm components that can affect the overall rate of convergence.  We have seen cases
where a structure preserving method might outperform an auxiliary method or vice versa.
Though both the auxiliary and structure preserving methods are generally successful algebraic multigrid (AMG) solvers, there are still important robustness issues for complex problems with large material variation
and stretched meshes.
This includes cases where the convergence of existing solvers might stagnate before reaching an acceptable accuracy. In other scenarios the solver might converge but
the number of iterations increases noticeably as one refines the mesh or more generally the solver requires many more iterations than one would expect from a multigrid solver.

In this paper, we propose a new AMG solver that can also leverage standard AMG solvers designed for nodal Poisson-like problems.
Like the auxiliary methods, the new solver can easily incorporate weak connection dropping, which is very important for anisotropic or material-variation problems.
As we will show, the new solver has close connections to geometric multigrid and does not require a second AMG sweep within the preconditioner
to correct for deficiencies.

%%% Local Variables:
%%% mode: LaTeX
%%% TeX-master: "hcurlEmin"
%%% End:

%%%%%%%%%%%%%%%%%%%%%%%%%%%%%%%%%%%%%%%%%%%%%%%%%%%%%%%%%%%%%%%%%%%%%
%%%%%%%%%%%%%%%%%%%%%%%%%%%%%%%%%%%%%%%%%%%%%%%%%%%%%%%%%%%%%%%%%%%%%
%%%%%%%%%%%%%%%%%%%%%%%%%%%%%%%%%%%%%%%%%%%%%%%%%%%%%%%%%%%%%%%%%%%%%
\section{Reitzinger and Sch\"oberl}\label{sec:rsAMG}
\begin{figure}[htb!]\label{fig:rs example}
\centering
\vspace*{-0.1in}
%{\includegraphics[trim=200 250 400 50,clip,scale=.4]{Figs/PiecewiseConst2.png}}
 {\includegraphics[trim=200 250 400 50,clip,scale=.4]{     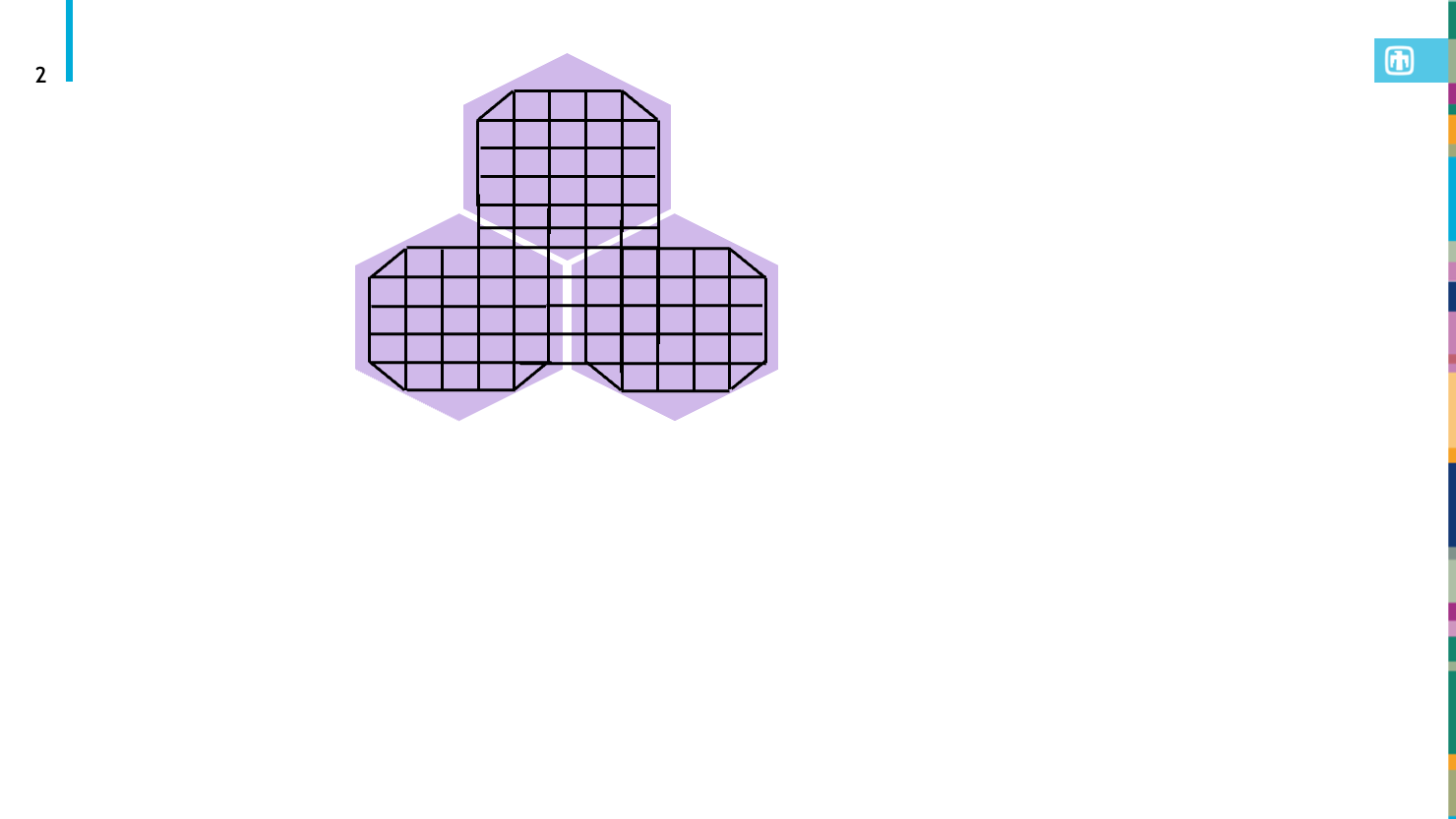}}
%{\includegraphics[trim=370 200 300 150,clip,scale=.6]{Figs/PiecewiseConstRS.png}}
 {\includegraphics[trim=370 200 300 150,clip,scale=.6]{     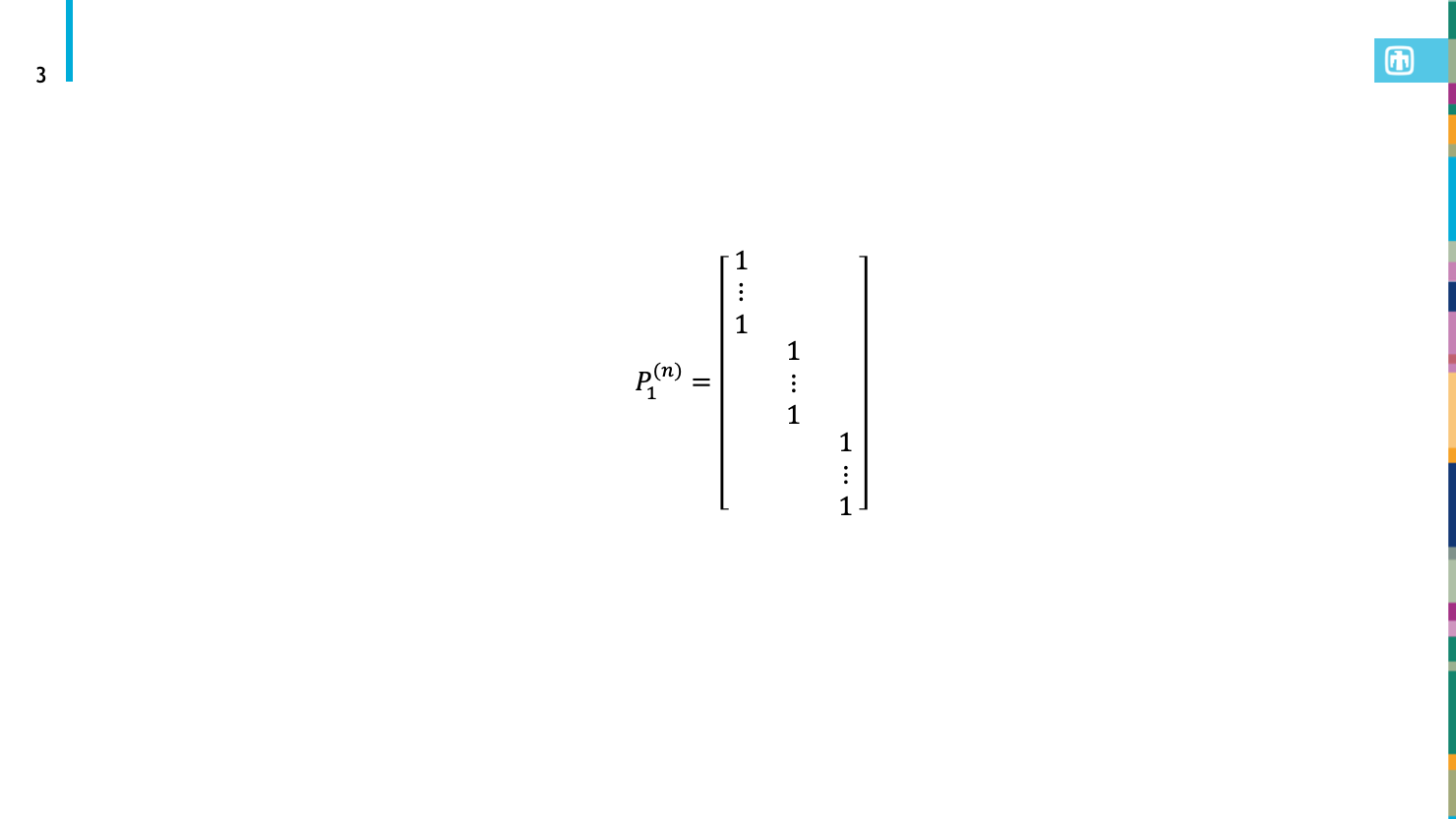}}
\caption{Piece-wise constant interpolation example. Purple hexagons illustrate three aggregates each corresponding to one column of
$P^{(n)}_1$. The $1^{st}$ and $2^{nd}$ columns each consist of 33 non-overlapping nonzeros while the $3^{rd}$ column contains 34 nonzeros.}
\end{figure}

The Reitzinger and Sch\"oberl AMG algorithm (RSAMG) begins by first
considering a finite element discretization of
a related nodal scalar PDE given by
\begin{equation} \label{eq: nodal problem}
-\Delta u  + \sigma(x,y,z) ~  u = f
\end{equation}
and generates an AMG hierarchy for this nodal problem. In what follows, we describe a two-level hierarchy to simplify notation.
However, the ideas are easily extended to multi-level hierarchies. The multi-level AMG hierarchy is effectively
described by a set of nodal interpolation operators $P^{(n)}_\ell$. For a two-level hierarchy, there is just
one operator that interpolates corrections from the coarse grid
to the fine grid and so the subscript $\ell$ is dropped from the description.
This interpolation operator correspond to non-overlapping piece-wise constant
interpolation. That is, each row of $P^{(n)}$ contains exactly one nonzero
and the value of this nonzero is always one. Normally, these piece-wise
constant grid transfers are constructed using an aggregation process such as
that used for smoothed aggregation~\cite{Vanek1996,Vanek2001}. Figure~\ref{fig:rs example} illustrates a mesh
with a set of aggregates and an associated piece-wise constant grid transfer.
Effectively, one can view each aggregate as a node on the coarse grid.

To build a coarse mesh, coarse edges must be defined that effectively connect the coarse nodes.
The Reitzinger and Sch\"oberl procedure for constructing edges is equivalent to first
projecting a Poisson operator
\begin{equation} \label{eq: RS coarse gradient}
Z = [ P^{(n)}]^T ~~ [ D^{(n\rightarrow e)}_h]^T [ D^{(n\rightarrow e)}_h]~~ P^{(n)} ,
\end{equation}
and then defining one coarse edge (between coarse nodes $i$ and $j$) for each
off diagonal nonzero $Z_{ij}$ in the strictly upper triangular portion of $Z$.
Notice that the projected Poisson operator is defined using the discrete gradient $D^{(n\rightarrow e)}_h$.
By careful examination, one can see that
the $j^{th}$ column of $ D^{(n\rightarrow e)}_h P^{(n)}$
corresponds to interpolating the $j^{th}$ coarse vertex
to the $j^{th}$ aggregate via piece-wise constant interpolation followed by applying a discrete
gradient operator. The resulting vector has nonzero entries only along fine edges
with exactly one fine vertex in the $j^{th}$ aggregate. As $Z_{ij}$ is defined by taking
inner-products between the $i^{th}$ and $j^{th}$ columns of
$ D^{(n\rightarrow e)}_h P^{(n)}$, it
follows that
the strictly upper triangular nonzeros of $Z$ give rise to a set of coarse edges that
connect adjacent aggregates.
In this sense,
the algorithm has a  relationship with a Delaunay tessellation which creates edges by joining adjacent
Voronoi regions. While aggregates are not Voronoi regions, we find this connection interesting and
it does give some sense that if the aggregates are fairly regularly shaped, then a reasonable
set of coarse edges will be defined.  Notice that in Figure~\ref{fig:rs example}, three coarse edges
are created connecting adjacent aggregates and if all aggregates happen to be hexagonal in a larger mesh,
a reasonable set of coarse edges would be defined by this procedure.

Overall,
\eqref{eq: RS coarse gradient} allows one to construct
a coarse $D^{(n\rightarrow e)}_H $
matrix using only a nodal interpolation operator and a user-provided
$D^{(n\rightarrow e)}_h$. In the multi-level case, we would generate a hierarchy
of discrete gradient operators  using a hierarchy of nodal interpolation operators.
Before discussing edge interpolation within RSAMG, we note that special treatment is
needed to address Dirichlet boundary conditions. Dirichlet boundary conditions give rise to $D^{(n\rightarrow e)}_h $ edges
with only one nonzero corresponding to the one interior node so that $ S^{(e)}_{H} D^{(n\rightarrow e)}_h = 0$
even when Dirichlet boundary conditions are present. These single-node edges must also be projected to the
coarse mesh when defining $D^{(n\rightarrow e)}_H $. This is done by simply taking each coarse node that interpolates to
the interior node of a single-node edge and forming a single-node coarse edge (one edge for each coarse node) and then augmenting
$D^{(n\rightarrow e)}_H $ with these additional coarse edges.

To construct an AMG hierarchy for the original matrix problem
\eqref{eq: discrete eddy}, Reitzinger and Sch\"oberl define an edge interpolation operator
$P^{(e)}$ such that the following commuting relationship holds
\begin{equation} \label{eq: commute}
P^{(e)} D^{(n\rightarrow e)}_{H} =
D^{(n\rightarrow e)}_h P^{(n)}  .
\end{equation}
We use \(P^{(e)}\) to define a coarse grid operator \(A^{(e)}_{H}\) via Galerkin projection
$$
A^{(e)}_H =  [P^{(e)}]^T ~~A^{(e)}_h ~~ P^{(e)} .
$$
This coarsening is structure preserving in the sense that if we similarly define
$$
S^{(e)}_H =  [P^{(e)}]^T ~~S^{(e)}_h ~~ P^{(e)}
$$
as the coarse-grid approximation of the curl-curl term, then \eqref{eq: commute} implies that
\begin{align}
S^{(e)}_H  D^{(n\rightarrow e)}_{H} &=    [P^{(e)}]^T S^{(e)}_h P^{(e)} D^{(n\rightarrow e)}_{H} \\
                                               &=    [P^{(e)}]^T S^{(e)}_h D^{(n\rightarrow e)}_h P^{(n)}  \nonumber \\
                                               &=    \Theta_H . \nonumber
\end{align}
Hence the property \eqref{eq: discrete null space condition} satisfied by \(S^{(e)}_h\) is preserved.
Here, $\Theta_H \in \mathbb{R}^{n_H \times m_H}$
is a zero matrix  where $n_H$ is the number of coarse edges while $m_H$ is the number of coarse vertices.
The second line follows from the commuting relationship while the third line follows by assuming that %from an inductive argument where we assume that
on the fine level we have $S^{(e)}_h D^{(n\rightarrow e)}_h = \Theta$.
As \cite{ReSc02} show, the $P^{(e)} $ that satisfies this relationship has columns (or basis functions) that are only nonzero (with values of
$\pm 1.0$) corresponding to fine edges that connect adjacent aggregates, which is related to the nonzero structure of
$ D^{(n\rightarrow e)}_h P^{(n)}$ discussed in the description of $Z$ above.
We omit the details and refer interested readers to \cite{ReSc02}. The main drawback of RSAMG is its reliance on sub-optimal
piecewise-constant nodal interpolation which leads to sub-optimal edge interpolation where coarse edges only interpolate to
fine edges that span two aggregates. Thus, many rows of $P^{(e)}$ have no nonzeros.

%%%%%%%%%%%%%%%%%%%%%%%%%%%%%%%%%%%%%%%%%%%%%%%%%%%%%%%%%%%%%%%%%%%%%
%%%%%%%%%%%%%%%%%%%%%%%%%%%%%%%%%%%%%%%%%%%%%%%%%%%%%%%%%%%%%%%%%%%%%
%%%%%%%%%%%%%%%%%%%%%%%%%%%%%%%%%%%%%%%%%%%%%%%%%%%%%%%%%%%%%%%%%%%%%
\section{Energy minimization AMG background}\label{sec:eAMG_background}
The H-curl AMG method described in Section~\ref{sec:eminPlusHcurl} leverages an energy-minimizing framework for generating interpolation operators
that was proposed in~\cite{Olson2011} and further developed in~\cite{JaFrScOl2023}. It is a generalization of ideas
in~\cite{Mandel1999} and related to similar schemes in
\cite{Brandt00generalhighly,Brandt01multiscalescientific, BrZi2007, KoVa2006, Va2010, Wa2000, WaChSm2000,XuZi2004}.
We briefly summarize the most relevant aspects of this framework for the H-curl AMG method that will be presented in Section~\ref{sec:eminPlusHcurl}

Let ${\cal N}$ be a set of
$m_h \times m_H$
matrices with a specific nonzero pattern,
$M_h$ be a set of fine level modes requiring exact
interpolation, and $M_H$ be a coarse representation of these modes (e.g., obtained by injection of the fine modes onto the coarse nodes).
Prolongator coefficients are determined through an
approximate solution to a constrained quadratic minimization problem
\begin{equation}
%  P = \argmin\limits_P \sum\limits_j \frac{1}{2} \|P_j\|_{A}^2
  P = \argmin\limits_P \sum\limits_j \frac{1}{2} \|P_{:j}\|_{A}^2
  \label{eq:min_problem} \\
%\end{equation}
~~~~~~\mbox{subject to}~ \\
%\begin{equation}
  P \in {\cal N},\quad\mbox{and}\quad M_h = P M_H .
\end{equation}
Here, $A$ is the matrix system that one wishes to solve via AMG while $P_{:j}$ refers to the $j$th column of $P$, and the sum is over all columns in $P$.
%In this section, we only discuss a two-level method and do not use subscripts on $A$ nor $P$ to denote the level within the hierarchy.
The modes $M_h$ and $M_H$ are typically a collection of smooth vectors that are problematic for standard iterative solvers such as
the constant vector for Poisson problems or rigid body modes for linear elasticity. These could also be modes that have been computationally
identified by some iterative process as within an adaptive AMG scheme, e.g.~\cite{BrFaMaMaMcRu2005,bootstrap,NMTMA-8-1,10.1145/3190647}.
Overall, the objective function favors smooth basis functions while the constraints assure the modes problematic for fine level relaxation
are represented on the coarse grid.

This quadratic minimization problem is equivalent to solving the linear system
\begin{equation} \label{eq: KKT}
\left [
\begin{matrix}
\hat{A}  & X^T \\
X  & 0
\end{matrix}
\right ]~\left [
\begin{matrix}
             \widebreve{\delta P} \\
\lambda
\end{matrix}
\right ]
=
\left [
\begin{matrix}
-\widebreve{AP_0} \\
0
\end{matrix}
\right ]
\end{equation}
where $\hat{A}$ is a block matrix with $A$ repeated $m_H$ times along the block diagonal, $X$ is a matrix representing the two types of linear constraints, $P_0$ is a feasible initial guess satisfying
the constraints, and $\delta P$ is a correction to $P_0$ giving the desired optimal solution $P = P_0 + \delta P$.
%and finally $P = P_0 + \delta P$.
 The operation $y =             \widebreve{Y} $ converts the matrix ${Y}$ to a vector by
inserting all matrix entries ${Y}_{ij} $ into $y$ at the
${(j\hskip -.03in
-\hskip -.03in 1)\hskip -.03in*\hskip -.03inm_h\hskip -.03in+\hskip -.03ini}$ location.
A suitable $P^{(0)}$ can be obtained by solving an under-determined linear system involving $X$ as discussed in ~\cite{Olson2011}.
The saddle-point system~\eqref{eq: KKT} can be solved by the projected conjugate gradient method which effectively applies CG to the matrix
$(I - X^T (X X^T)^{\dagger} X ) \hat{A} (I - X^T (X X^T)^{\dagger} X )$ and has the added benefit of minimizing the error in the $A$-norm.
Here, $(X X^T)^{\dagger}$ denotes the pseudo-inverse of $X X^T$ to address cases when $X X^T$ is singular.
In an implementation, there is no need to explicitly form $\hat{A}$ or explicitly apply a matrix-to-vector conversion by instead leveraging
matrix-matrix products. Further,
the second projection can be omitted as %it is equivalent to applying the identity as
 the previous approximate solution to $            \widebreve{\delta P}  $ satisfies the constraint
when the starting initial guess for $\delta P$ is zero.

The smoothed aggregation AMG method~\cite{Vanek1996,Vanek2001} can be viewed as applying a one-step damped Jacobi iteration
to solve~\eqref{eq: KKT} using a tentative prolongator as $P_0$.  That is,
\begin{equation} \label{emin iteration}
 \delta P = D^{-1} A P^{(k)} ; \hskip .1in \delta p = (I - X^T (X X^T)^{\dagger} X )             \widebreve{\delta P}  ; \hskip .1in
P^{(k+1)} \leftarrow P^{(k)} - \omega [\overset{     \frown}{\delta p}]
\end{equation}
where
$D$ is the diagonal of $A$; $P^{(k)}$ refers to the $k^{th}$ approximation to $P$.
 The vector-to-matrix conversion $[\overset{\frown}{y}] $ is the reverse operation of $            \widebreve{Y} $.
The damping parameter $\omega$ is usually chosen to be $\frac{4}{3} \hat{\rho}(D^{-1} A)$
where $\hat{\rho}(D^{-1} A)$ is an approximation to the maximum eigenvalue of $D^{-1} A$.
However in the smoothed aggregation context,
the constraints are automatically or nearly automatically satisfied due to the choice of $P^{(0)}$,
and so the projection step is omitted. However, the projection must be retained for a more general energy minimization algorithm.

While solving a constrained optimization problem could be computationally expensive compared to the cost of solving a linear system,
the solution to this optimization problem can be very roughly approximated. For example, smoothed aggregation takes only one iteration
and more generally we find that at most only 1-4 iterations are needed to generate
suitable grid transfers leading to effective AMG convergence rates. Additionally, a prolongator from a previous linear solve
can often be used as a good initial guess to the energy minimization process when solving a sequence of related problems
that might occur within a nonlinear solver or when implicit time stepping is employed.
The most expensive operation within a single iteration is applying the projection
operation due to the need to perform a QR factorization involving $X$ or invert
the matrix $X X^T$.  While this might seem prohibitive, the projection for sparsity pattern constraints
corresponds to simply setting any nonzeros outside of ${\cal N}$ to 0.
For mode constraints, $X$ has a block diagonal structure due to the fact that constraints
for one row of $\delta P$ are self-contained or completely independent of constraints for
any other row of $\delta P$.  Thus, several small independent QR factorizations or matrix inverses
need to be performed as opposed to a large globally coupled factorization. The row dimension
of the blocks is equal to the column dimension of $M_h$.  The column dimension of a block
used to update the $i^{th}$ row of $\delta P$ is equal to the number of nonzeros in the associated row of $\cal N$.
Thus, if the number of modes is modest (e.g., one for Poisson operators or six for linear elasticity)
and the number of nonzeros per row in $\cal N$ is not large, enforcement of the constraints is not prohibitively
expensive.  This is discussed in~\cite{JaFrScOl2023} which includes an evaluation of the overall setup cost
along with several other practical aspects
of energy minimization AMG.
A precise example of $X$ is provided in Section 2.1 of ~\cite{Olson2011} with some
linear algebra details.

\section{An Energy minimization Algorithm for H-curl}\label{sec:eminPlusHcurl}
From Section~\ref{sec:eAMG_background}, the energy-minimization framework
seeks to approximately solve
\begin{equation}
  P = \argmin\limits_P \sum\limits_j \frac{1}{2} \|P_{:j}\|_{A}^2
  \label{eq:min_hcurl} \\
~~~~~~\mbox{subject to}~ \\
  P \in {\cal N},\quad\mbox{and}\quad M_h = P~ M_H .
\end{equation}
By taking
$M_h = D^{(n\rightarrow e)}_h P^{(n)} $
and
$M_H = D^{(n\rightarrow e)}_{H}$.
This can be re-written as
\begin{equation}
  {P^{(e)}} \hskip -.05in = \hskip -.03in \frac{1}{2} \argmin\limits_{P^{(e)}} \hskip -.03in \sum\limits_{j} \|P^{(e)}_{:j}\|_{S^{(e)}_h}^2
  \label{eq:min2_hcurl} \\
~\mbox{subject to}~ \\
  {P^{(e)}} \hskip -.06in \in {\cal N},\mbox{ and } D^{(n\rightarrow e)}_h P^{(n)} \hskip -.04in = \hskip -.04in P^{(e)} D^{(n\rightarrow e)}_{H}
\end{equation}
where we have also substituted $S^{(e)}_h$\footnote[4]{$S^{(e)}_h$ is semi-positive definite and so it does not define a norm, though
this is a minor technicality that does not affect practical application of an iterative approximation to
$P^{(e)}$. One can consider instead using $A^{(e)}_h$  or alternatively adding a gauge term to $S^{(e)}_h$ to avoid this mathematical nuance.} for $A$, and $P^{(e)}$ for $P$.
That is, the energy-minimization framework can be easily adapted to minimize the energy of interpolation basis functions subject
to the commuting constraints. The key point is that this can be accomplished for more general nodal interpolation
operators than the sub-optimal piece-wise constant interpolation considered by Reitzinger and Sch\"oberl. Of course, this requires
application of the projection step discussed in the previous Section.
The main differences with the general energy-minimization framework is that the operator $D^{(n\rightarrow e)}_{H}$ must be
defined as input for the energy-minimization procedure and that now the constraint equations have a matrix form as opposed to vectors
or modes.  % that must be accurately interpolated.
Normally, the column dimension of $D^{(n\rightarrow e)}_{H}$ will be considerably larger than a standard $M_h$ defined
by a few modes. Because of this, it may seem that the energy-minimization projection step might be prohibitively
expensive.  However, another key difference is that a standard $M_h$ is dense while $D^{(n\rightarrow e)}_{H}$ is very sparse
(at most two nonzeros per row with entries that are either one or minus one) and this will help reduce
the associated costs of the projection step. One final difference between the standard energy minimization
case and our specific context concerns rank/conditioning issues associated with the energy minimization 
projection step, $I - X^T (X X^T)^{\dagger} $X. As the constraints are now governed by $D^{(n\rightarrow e)}_{H}$,
we will see that there is a directed graph interpretation of $X$ and that $X X^T$ is essentially a graph Laplacian
whose properties can be leveraged.  Specifically,
there is no need to employ a singular value decomposition or rank revealing QR algorithm
to address potentially ill-conditioning or singularity issues, which might be needed in other energy minimization scenarios (e.g., for elasticity).
To complete the description of a structure preserving
energy-minimization algorithm applied to H-curl problems (SpHcurlAMG), we must describe
an algorithm for computing $D^{(n\rightarrow e)}_{H}$ and the sparsity pattern of
$P^{(e)}$. Further, we must look more closely at the constraint equations in terms of feasible solutions
and investigate practical computational aspects toward computing the projection.

\subsection{Computation of $D^{(n\rightarrow e)}_{H}$}
As with the RSAMG algorithm, the right hand side of the commuting equation defines a set of
coarse nodes that must be connected with coarse edges that will effectively define $D^{(n\rightarrow e)}_{H}$.
We recall that the Reitzinger and Sch\"oberl method for constructing $D^{(n\rightarrow e)}_{H}$
involved a projection of the form
$$
[ P^{(n)}]^T ~~ [ D^{(n\rightarrow e)}_h]^T [ D^{(n\rightarrow e)}_h]~~ P^{(n)}
$$
which has a nice analogy with Delaunay tessellations when $P^{(n)}$ is based on piecewise-constant interpolation.
Unfortunately, when $P^{(n)}$ is a more sophisticated nodal interpolation operator, there is no analogy with
Delaunay tessellation and in general we find that such an algorithm generates too many coarse edges where some of these
edges might also cross each other. If coordinates are supplied, one could apply a Delaunay algorithm to create a
tetrahedral mesh (somehow discarding tetrahedra that lie outside of a concave physical domain).

When a nodal prolongator is constructed by employing
smoothed aggregation (SA), a tentative piece-wise constant prolongator $ \bar{P}^{(n)}_\ell$ is first constructed
and then improved via a prolongator smoothing step. Another possibility for creating
$D^{(n\rightarrow e)}_{H}$ would be to retain $ \bar{P}^{(n)}_\ell$ and employ it
in conjunction with the Reitzinger and Sch\"oberl algorithm for creating coarse edges.
That is, form
\begin{equation} \label{eq: revised RS DH}
[ \bar{P}^{(n)}_\ell]^T ~~ [ D^{(n\rightarrow e)}_h]^T [ D^{(n\rightarrow e)}_h]~~ \bar{P}^{(n)}_\ell
\end{equation}
to then define $D^{(n\rightarrow e)}_{H}$ while still using the more sophisticated smoothed prolongator
${P}^{(n)}_\ell$ in the energy-minimization constraints. In this way, the Delaunay analogy still holds,
which should produce more physically realistic coarse meshes.
By developing an algorithm that converts a ${P}^{(n)}_\ell$ to a piece-wise constant $\bar{P}^{(n)}_\ell$,
one can generalize the basic idea to allow the possibility of using other AMG algorithms (besides smoothed
aggregation) for the purposes of constructing $P^{(n)}$. Such  a conversion method 
is illustrated in
Algorithm~\ref{algo: constant interpolation} where the basic idea is to define a sparsity pattern
\begin{algorithm}
\begin{tabbing}
\hskip .4in \=  \hskip .17in \=  \hskip .17in \=  \hskip .17in \=  \kill  \\ \linefill \\
$P_{const}$ = \textsc{Gen\_Piecewise\_Const}($P$, nPasses) \\[2pt]
~~~~~ ! Construct $ P_{const} \in \mathbb{R}^{m_h \times m_H}$ by assigning each row to one column. \\[2pt]
~~~~~ ! Make multiple sweeps assigning a few largest magnitude entries. Safeguards\\
~~~~~ ! (not shown) needed  to guarantee each column has some row assignments. \\[10pt]

\> targetNzPerCol $\leftarrow$ Choose\_Target\_Proportional\_To\_Ps\_NnzPerCol(P) \\[2pt]

\> col\_Assignment$(1:m_h) \leftarrow 0$   \\[7pt]

\> for $k=1:\mbox{nPasses}$ \{ \\[2pt]
\>\> for $j=1:m_H$ \{ \\[2pt]
\>\>\> nAssign $\leftarrow \left \lceil{\mbox{targetNzPerCol}/\mbox{nPasses}} \right \rceil $ \\[2pt]
\>\>\> rows $\leftarrow $ Largest\_Unassigned($P_{:j}$,col\_Assignment, nAssign ) \\[2pt]
\>\>\> col\_Assignment(rows) $\leftarrow j$ \\[-2pt]
\>\> \} \\[-2pt]
\> \}\\[7pt]

\> ! employ greedy assignment to any unassigned rows \\[5pt]

\> col\_Assignment $\leftarrow$ Greedy\_Assign(col\_Assignment) \\[4pt]

\> for $i=1:m_H$ \{ $[P_{const}]_{i,\mbox{col\_Assignment(i)}} \leftarrow 1 $ \}
\end{tabbing}
\caption{\label{algo: constant interpolation} Pseudo-code to construct piece-wise constant interpolation}
\end{algorithm}
for $\bar{P}^{(n)}_\ell$ corresponding to non-overlapping basis functions by dropping small entries,
allowing one to consider
classical AMG schemes~\cite{BrMcRu1984,RuSt1987} to generate $P^{(n)}$.
Even though we employ smoothed aggregation to generate the $P^{(n)}$, we use
Algorithm~\ref{algo: constant interpolation} to construct $\bar{P}^{(n)}_\ell$ for the purposes of
the $D^{(n\rightarrow e)}_{H}$ computation.
There is one caveat that will be discussed shortly in the next sub-section.

\subsection{A closer look at the commutator constraints}
Recall the commutator equation
$$
P^{(e)} D^{(n\rightarrow e)}_{H} = D^{(n\rightarrow e)}_h P^{(n)}
$$
where it is assumed that the right hand side and $D^{(n\rightarrow e)}_{H}$ are defined
before initiating the energy minimization iteration.
Notice that $P^{(e)} D^{(n\rightarrow e)}_{H} v = 0 $ when $v$ is a constant
vector and all rows of $D^{(n\rightarrow e)}_{H}$ contain precisely two nonzeros. This is due to the
nature of the discrete gradient. This implies that $P^{(n)} $ must preserve constants
so that the right hand side of the commutator equation is also zero. While $D^{(n\rightarrow e)}_{H}$
may contain a few rows with only one nonzero (due to Dirichlet boundary conditions), we always modify
the nodal prolongator by scaling the coefficients in each row so that all row sums are identically equal
to one.

\begin{figure}[htb!]\label{fig:many pictures}
\begin{minipage}{.605\textwidth}
\centering
%{\includegraphics[trim=75 400 470 0,clip,scale=.4]{Figs/Ppattern1.png}}
 {\includegraphics[trim=75 400 470 0,clip,scale=.4]{     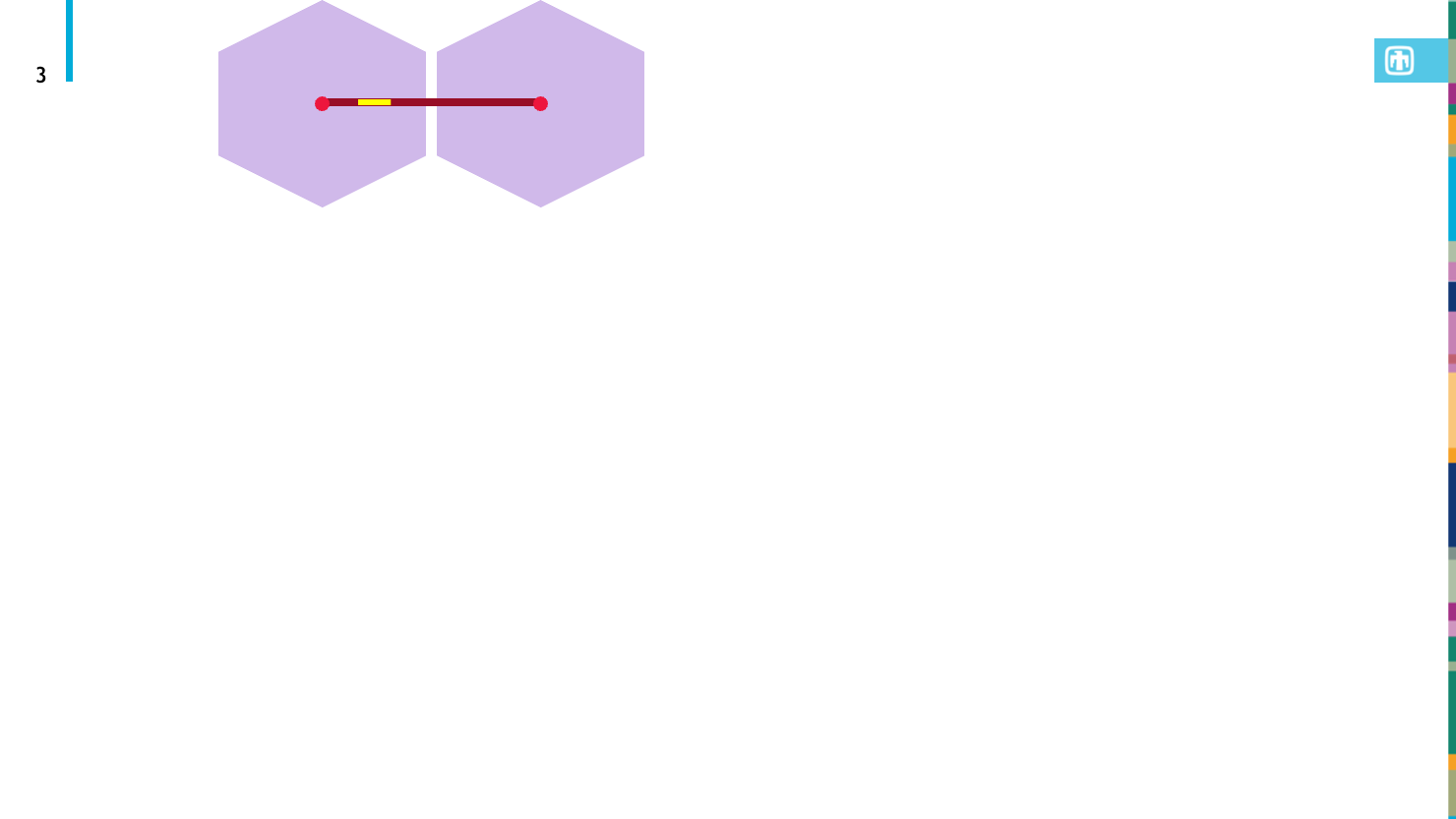}}
\end{minipage}
\begin{minipage}{.38\textwidth}
$\left [\begin{matrix}
1  & -1 % &     &     &
\end{matrix}  \right ]
$
\end{minipage}
\vskip .3in
\begin{minipage}{.6\textwidth}
\centering
%{\includegraphics[trim=75 300 470 0,clip,scale=.4]{Figs/Ppattern2.png}}
 {\includegraphics[trim=75 300 470 0,clip,scale=.4]{     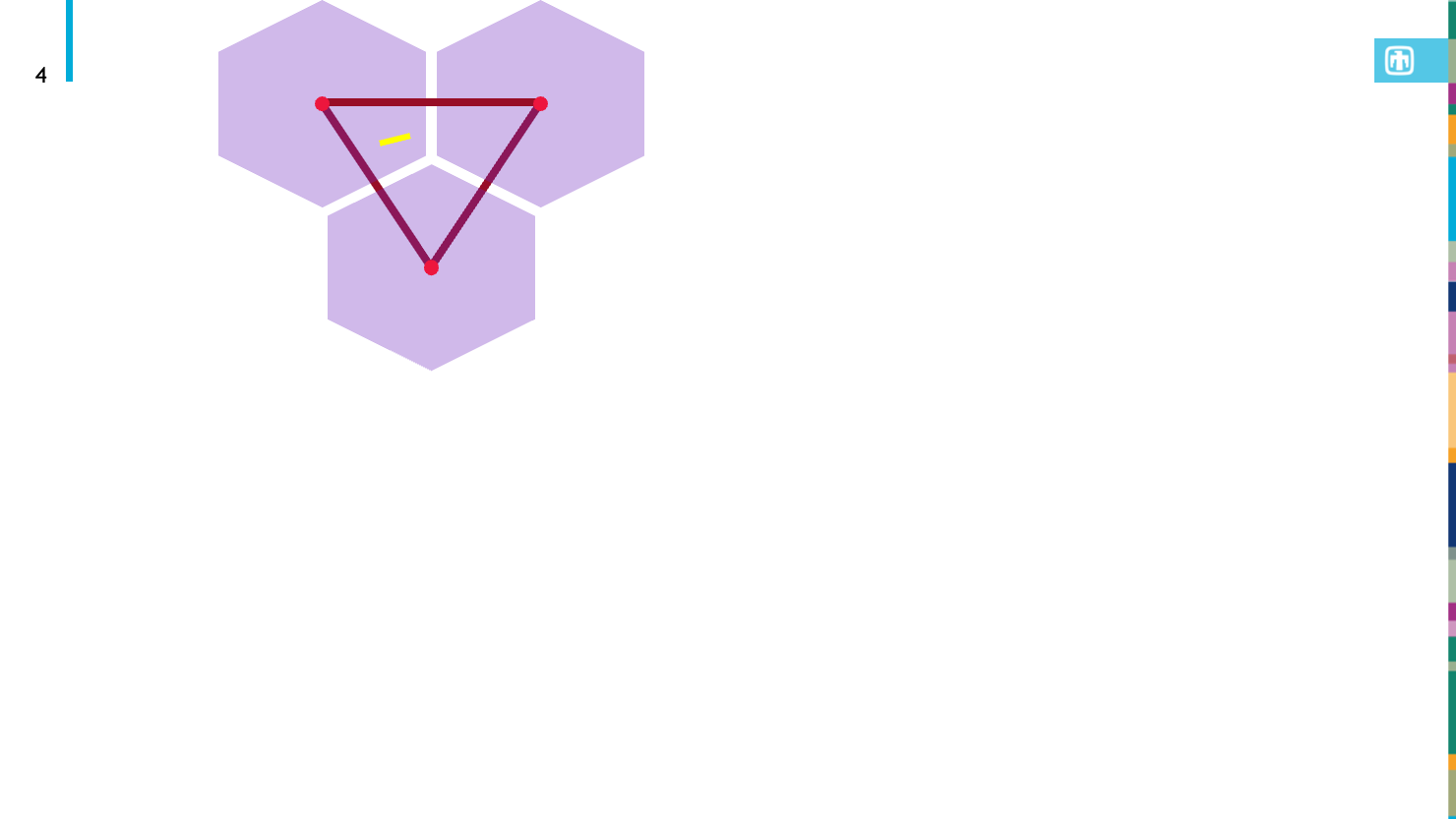}}
\end{minipage}
\begin{minipage}{.38\textwidth}
$\left [\begin{matrix}
\phantom{-}1  &             -1 &               \\
              &   \phantom{-}1 &           -1  \\
          -1  &                & \phantom{-}1
\end{matrix}  \right ]
$
\end{minipage}
\vskip .3in
\begin{minipage}{.6\textwidth}
\centering
%{\includegraphics[trim=145 300 275 0,clip,scale=.4]{Figs/Ppattern3.png}}
 {\includegraphics[trim=145 300 275 0,clip,scale=.4]{     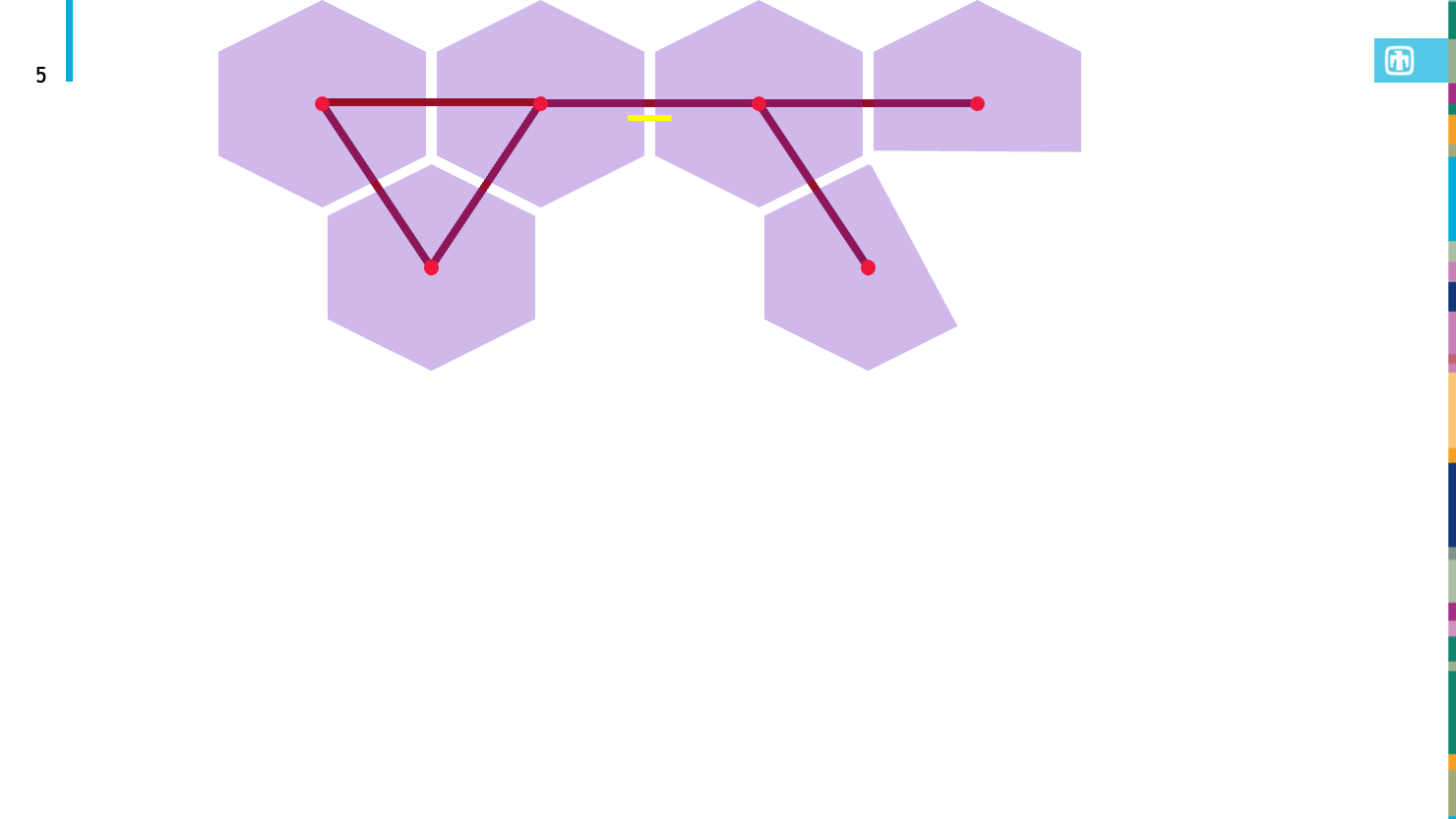}}
\end{minipage}
\begin{minipage}{.38\textwidth}
$\left [\begin{matrix}
%    1                2              3               4              5         6
\phantom{-}1  &             -1 &              &               &           &          \\
              &   \phantom{-}1 &           -1 &               &           &          \\
          -1  &                & \phantom{-}1 &               &           &          \\
              &                & \phantom{-}1 &           -1  &           &          \\
              &                &              & \phantom{-}1  &     -1    &          \\
              &                &              & \phantom{-}1  &           &    -1
\end{matrix}  \right ]
$
\end{minipage}
\vskip .3in
\begin{minipage}{.6\textwidth}
\centering
%{\includegraphics[trim=145 300 275 0,clip,scale=.4]{Figs/Ppattern5.png}}
 {\includegraphics[trim=145 300 275 0,clip,scale=.4]{     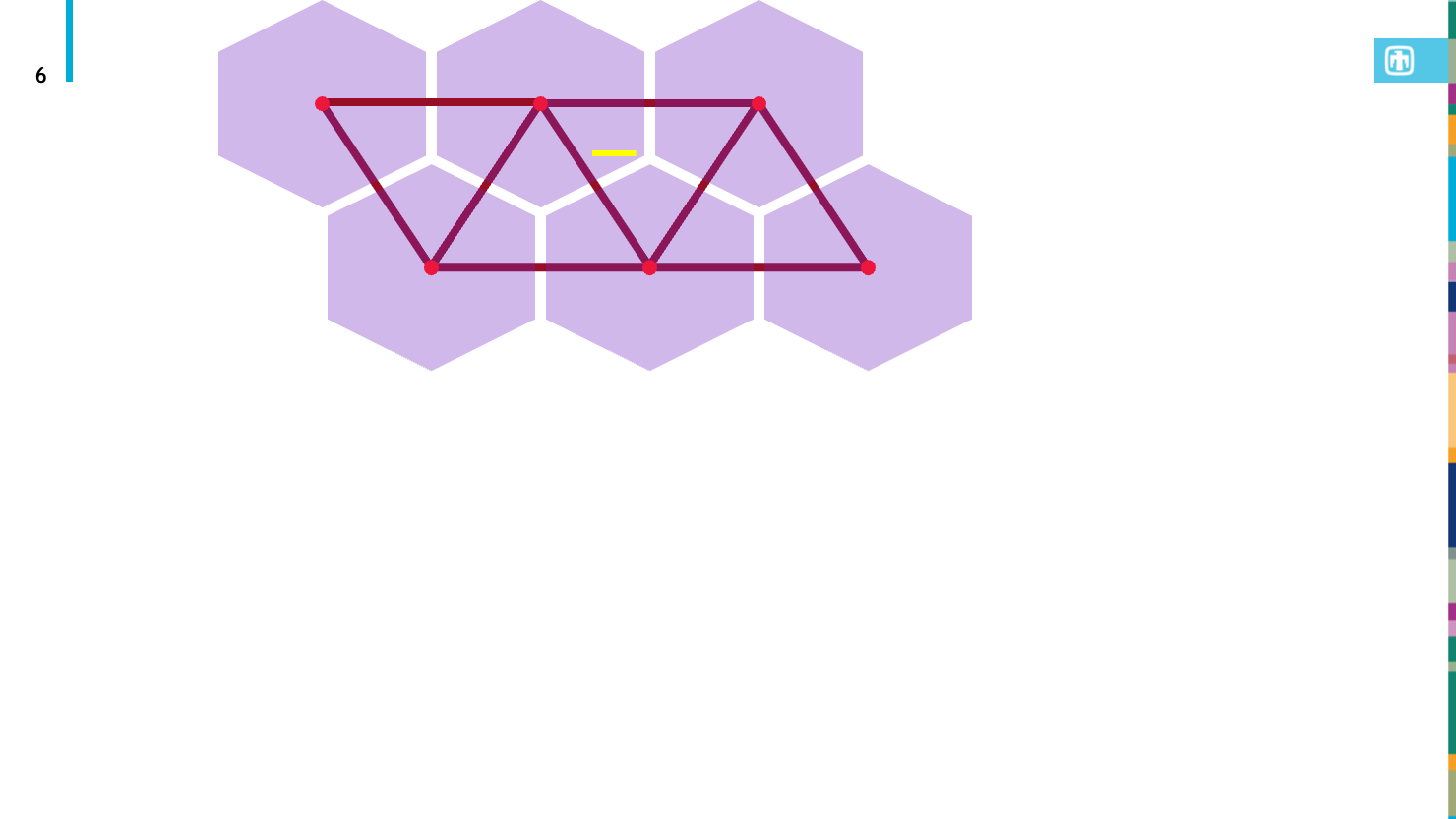}}
\end{minipage}
\begin{minipage}{.38\textwidth}
$\left [\begin{matrix}
%%       1             2              3               4               5         6
\phantom{-}1  &             -1 &              &               &              &              \\   % 1
              &   \phantom{-}1 &           -1 &               &              &              \\   % 2
          -1  &                & \phantom{-}1 &               &              &              \\   % 3
              &                &              & \phantom{-}1  &           -1 &              \\   % 4
              &                &              &               & \phantom{-}1 &           -1 \\  % 5
              &                &              &           -1  &              & \phantom{-}1 \\  % 6
              &   \phantom{-}1 &              &               &           -1 &              \\   % 2
              &                & \phantom{-}1 &           -1  &              &              \\   % 3
\end{matrix}  \right ]
$
\end{minipage}
\caption{Graphical illustration of prolongator sparsity pattern for four different scenarios.}
\end{figure}

We now examine the $i^{th}$ row of the commutator equation
written as
\begin{equation} \label{eq: one row}
p_{i:}~  D^{(n\rightarrow e)}_{H} = r_{i:}
\end{equation}
where
$ p_{i:} \in \mathbb{R}^{1 \times n_H}$,
$r_{i:} \in \mathbb{R}^{1 \times m_H}$ ,
$n_H$ is the number of coarse edges, $m_H$ is the number of coarse nodes,
and $r_{i:}$ is the $i^{th}$ row of $D^{(n\rightarrow e)}_h P^{(n)} $.
The sparsity pattern of $r_{i:}$ corresponds to all coarse nodes that interpolate to one (or perhaps both) of the two nodes
defining the $i^{th}$ edge\footnote{There is some chance
that a sparsity pattern entry might be {\it missing} due to cancellation (when a coarse node interpolates
to both edge end points with the same weight). We avoid this by adding a small random perturbation to either
$D^{(n\rightarrow e)}_h$ or $ P^{(n)} $ when constructing a sparsity pattern for $P^{(e)}$.}.
Figure~\ref{fig:many pictures} illustrates four different scenarios for the sparsity pattern of
$p_{i:}. $ In each leftmost image, purple denotes the support of the piece-wise constant version of the
nodal interpolation operator. The entire fine mesh is not shown to avoid clutter, but one should imagine
an underlying fine grid as was depicted in Figure~\ref{fig:rs example}. The $i^{th}$ fine grid edge
is, however, depicted in yellow. The red dots denote the nonzeros or coarse nodes in $r_{i:}$.
The dark red lines depict suitable coarse edges defining the sparsity pattern of $p_{i:}$. Thus,
the 4 images each
depict a fine edge that interpolates from $1, 3, 6$, and $9$ coarse edges respectively as one descends
from the top-most to the bottom-most image. The associated matrix just to the right of each image
corresponds to the relevant sub-matrix of $D^{(n\rightarrow e)}_{H} $. That is,
the sub-matrix rows correspond to nonzero coarse edges in $p_{i:}$ while the
sub-matrix columns correspond to the nonzero coarse nodes in $r_{i:}$.
In setting up the projection, a linear algebra problem involving this
sub-matrix must be computed. This can be seen by simply transposing ~\eqref{eq: one row}
to the form $[D^{(n\rightarrow e)}_{H}]^T  p_{i:}^T = r_{i:}^T $ so that the
desired unknown vector $ p_{i:}^T$ is applied on the right side of the matrix. The lowest image whose
transposed sub-matrix is $6 \times 9$  is closest to the most likely scenario.
The rank is $5$ and not $6$
due to the fact that the constant vector is in the null space of the discrete gradient sub-matrix, but importantly $r_{i:}^T $
will always be consistent so long at the nodal interpolation operator preserves constants.
As the system is under-determined, an infinite number of solutions exist.  A solution with
minimum 2-norm is given by $\tilde{Q}_i \tilde{R}_i^{-T} \tilde{r}_{i:}$ which can be found after applying a QR decomposition
to the $i^{th}$ $D^{(n\rightarrow e)}_{H} $ sub-matrix.
Here, $\tilde{Q}_i \in \mathbb{R}^{n_p \times k}$, $\tilde{R}_i \in \mathbb{R}^{k \times k}$, and $\tilde{r}_{i:}^T \in \mathbb{R}^{k \times 1}$
are sub-matrices of the QR factors or sub-vectors of ${r}_{i:}$ where $n_p$ is the number of nonzeros in $p_{i:}^T$'s sparsity
pattern and $k$ is the column rank of the $D^{(n\rightarrow e)}_{H}$ sub-matrix. Normally,
this column rank is one less than the column dimension of this sub-matrix as in our examples.
However, it may be equal to the the column dimension if {edges incident to Dirichlet nodes} are included
within the sub-matrix (discussed below).
The transposed sub-matrix for the two middle images are also under-determined even though
they are square because of the null space and consistent right hand side.
The top-most sub-matrix is $2 \times 1$.  This system can also be solved but
now has a unique solution (again due to the null space and consistent right hand side).
In this case, $p_{i:}$ has only one nonzero and so its value is effectively determined
only by the constraint and not the minimization process (as there is only one feasible
solution). One can verify that if the nodal interpolation operator happens to be given by linear interpolation in this one-coarse edge case,
then because of the constraints $P^{(e)}$ will effectively take the current coarse edge approximate solution
and scale it by the ratio of the lengths between the fine and coarse edge to define the interpolated fine
solution on the associated fine edge.

\begin{figure}[htb!]\label{fig: two pictures}
\begin{minipage}{.705\textwidth}
\centering
%{\includegraphics[trim=115 300 270 0,clip,scale=.4]{Figs/Ppattern4.png}}
 {\includegraphics[trim=115 300 270 0,clip,scale=.4]{     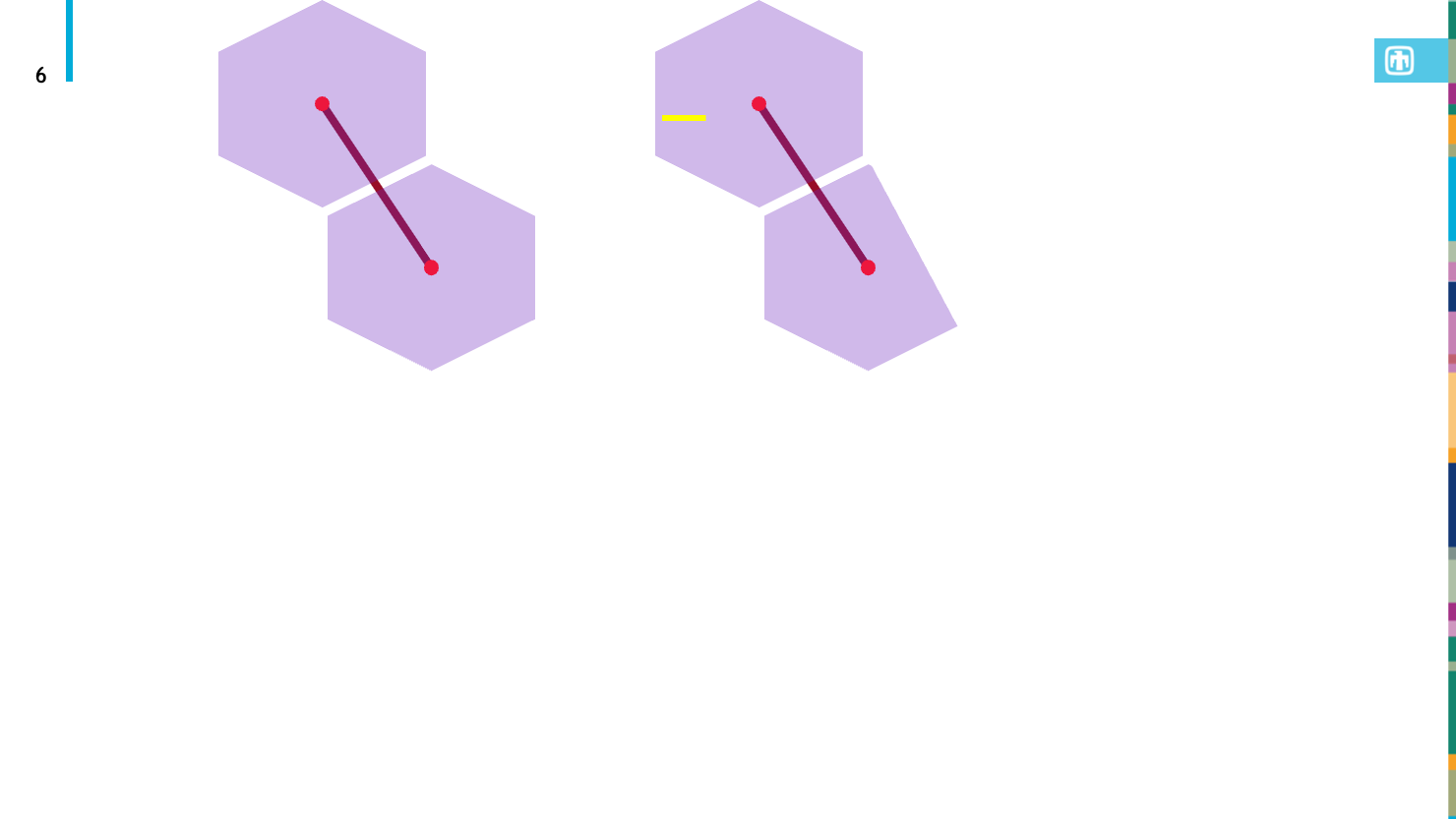}}
\end{minipage}
\begin{minipage}{.28\textwidth}
$\left [\begin{matrix}
           1  &  -1  &     &      \\
              &      &   1 &  -1
\end{matrix}  \right ]
$
\end{minipage}
\caption{Graphical illustration of prolongator sparsity pattern for irreducible sub-graph}
\end{figure}

All of the previous cases represent situations where the commutator equation
can be solved exactly. Notice that in each image, $p_{i:}$'s coarse edges are always
adjacent to two coarse nodes associated with $r_{i:}$ and that every coarse node
has at least one adjacent coarse edge in $p_{i:}$. These are in fact necessary conditions
for a feasible solution to exist.
Figure~\ref{fig: two pictures} depicts a scenario where a feasible solution generally does not exist.
In this case, the sub-matrix has a null-space with dimension two and $r_{i:}$ is not necessarily consistent.
If we insist that the coarse edge graph be irreducible (or fully connected), then this scenario can be
avoided. As irreducibility can be easily verified, the algorithm that builds the sparsity pattern for
the edge prolongator and defines the least squares systems additionally modifies the discrete coarse
gradient produced by the Reitzinger and Sch\"oberl algorithm adding an additional coarse edge if needed
so that the local coarse-edge graph is fully connected. In our experience, this happens very rarely,
e.g., far less than $1\%$ of the prolongator rows.
For the fully-connected coarse edge sub-graph, the null space is always at most one. The only
cases where the sub-matrix null space is empty occur when a Dirichlet edge (with only one
nonzero in the associated $D^{(n\rightarrow e)}_{H}$) is present ,
which can be easily verified.

Algorithm~\ref{algo: setup emin} gives a high-level view of the procedure to determine the edge interpolation
\begin{algorithm}
\begin{tabbing}
\hskip .4in \=  \hskip .17in \=  \hskip .17in \=  \hskip .17in \=  \hskip 1.4in \= \kill  \\ \linefill \\
$[~{\cal N}, Q, R, D^{(n\rightarrow e)}_{H}~ ]$ = \textsc{Emin\_Setup}($~D^{(n\rightarrow e)}_h,D^{(n\rightarrow e)}_{H},P^{(n)}~$)\\[2pt]
\> $\Tau \leftarrow$ Sparsity\_Pattern($~|D^{(n\rightarrow e)}_h |~  | P^{(n)} |~$)  \\[3pt]
\> $W \leftarrow$Diag($~\vv{\bf 3} - |D^{(n\rightarrow e)}_h |~ \overrightarrow{\bf 1}~) $   \>\>\>\> ! vector $\vv{\bf k}$ has entries equal to $k$, $W_{ii}=2$ if   \\[2pt]
\> \>\>\>\> ! $i^{th}$ edge has 1 endpoint. Otherwise, $W_{ii} = 1$ \\[2pt]
\> $B \leftarrow  | \Tau | ~|  D^{(n\rightarrow e)}_{H} |^T ~ W$                       \\[2pt]
%W ~|  D^{(n\rightarrow e)}_{H} | ~ | \Tau^T | $                                         \\[2pt]
\> ${\cal N} \leftarrow$ Initial\_Pattern(~$B$~)                            \>\>\>\> ! $ {\cal N}_{ij} = 1 ~\mbox{if}~ B_{ij} = 2 $ \\[2pt]

\> for $i=1:n^{(e)}_H$ \{ \>\>\>\> ! $ P^{(e)} \in \mathbb{R}^{n^{(e)}_h  \times n^{(e)}_H}$ \\[2pt]
\>\> ${\cal I} \leftarrow$ \{ $k \in {\cal N}_{ik} $~\} \>\>\>! Sparsity pattern of $p_{i:}$ in \eqref{eq: one row} \\[2pt]
\>\> ${\cal J} \leftarrow \Tau_{i:} $                \>\>\> ! Sparsity pattern of $r_{i:}$ in \eqref{eq: one row} \\[2pt]
\>\> isSingular $\leftarrow $ false \\[2pt]
\>\> DirichletEdges $\leftarrow$ \{ $k \in {\cal N}_{ik}  \mid W_{kk} \equiv 2 $ \} \\[2pt]
\>\> if (~$\mbox{DirichletEdges} \equiv \emptyset $~ ) \{  isSingular $\leftarrow $true \} \\[2pt]
\>\> ${\cal G}   \leftarrow [~D^{(n\rightarrow e)}_{H}~]_{{\cal I J}} $ \\[2pt]
\>\> if (~$\lnot$irreducible\_Graph($\phantom{|}{\cal G}\phantom{|}$)~) \{ \\[2pt]
\>\>\>   [~${\cal G}, D^{(n\rightarrow e)}_{H} ~] \leftarrow  $ make\_irreducible(~${\cal G}, {\cal J}$~) \\[2pt]
\>\> \} \\[2pt]
\>\> $[~\tilde{Q}_i , \tilde{R}_i~] \leftarrow $ QR(~${\cal G}$~)  \\[2pt]
\>\> if (~isSingular~) \{ $[~\tilde{Q}_i , \tilde{R}_i~] \leftarrow $ Reduce\_Dimension($~\tilde{Q}_i , \tilde{R}_i~$)  ~\} \\[2pt]
\> \}\\[7pt]

\end{tabbing}
\caption{\label{algo: setup emin} Pseudo-code to construct $\cal N$, possibly add edges to $D^{(n\rightarrow e)}_{H}$, and compute
$QR$ factors for least-squares sub-problems}
\end{algorithm}
sparsity patterns, augment $D^{(n\rightarrow e)}_{H}$ with additional edges to guarantee that a feasible
solution exists, % satisfying the constraints,
and sets up the least squares problem.
In the algorithm, $|~\dot~|$ indicates taking the absolute value, which is used instead of the random perturbations discussed
earlier to avoid cancellations. The sparsity pattern of $B$ is also the sparsity pattern of
$|D^{(n\rightarrow e)}_h |~ | P^{(n)}|~ | D^{(n\rightarrow e)}_{H} |^T $, which maps coarse edges to
fine edges. The matrix operators defining the $(i,j)$ entry of $B$ take the $j^{th}$ coarse edge and map it to its coarse node
end points. This is then interpolated to a set of fine nodes that is then mapped back to fine edges. In most case,  $B_{ij} = 2$ implies
that the two end points of the $j^{th}$ coarse edge affect the interpolated nodal solution at the end-points of the $i^{th}$ fine edge.
However, the
diagonal weight matrix $W$ scales coarse edges with one Dirichlet end point so that these one-vertex coarse edges are also
included in the sparsity pattern associated with the $i^{th}$ fine edge. If the sub-matrix graph is not irreducible, then
additional coarse edges are added to ${\cal G}$ and $D^{(n\rightarrow e)}_{H}$ so that the sub-graph is fully connected.
One possibility for choosing edges is to take the largest magnitude entry of $[ P^{(n)}]^T ~~ [ D^{(n\rightarrow e)}_h]^T  D^{(n\rightarrow e)}_h~~ P^{(n)}$ among a set of candidate coarse edges. We emphasize that here $P^{(n)}$ is not the piecewise constant nodal
operator used to construct the initial $D^{(n\rightarrow e)}_{H}$. It is also important to note that
there are typically only a few candidate coarse edges and so there is no need to explicitly form
$[ P^{(n)}]^T ~~ [ D^{(n\rightarrow e)}_h]^T  D^{(n\rightarrow e)}_h~~ P^{(n)}$ when only a few entries need to be
evaluated.

Finally, we note before closing this section that Algorithm~\ref{algo: setup emin} produces a {\it maximal} sparsity pattern for $P^{(e)}$. From the
perspective of solving the commutator equation, it is possible to prune some edges from this sparsity pattern (effectively reducing the row dimension of the sub-systems)
so long as the remaining graph is irreducible and all coarse nodes have at least one incident edge.
For example, if coordinates are available, one could remove coarse edges from the pattern of a fine edge row if that coarse is nearly orthogonal to the fine edge.

We summarize the steps to construct the edge prolongator \(P^{(e)}\) outlined in this section in Algorithm~\ref{algo: construction of edge prolongator}.

\begin{algorithm}
\begin{tabbing}
\hskip .4in \=  \hskip .17in \=  \hskip .17in \=  \hskip .17in \=  \hskip 1.4in \= \kill  \\ \linefill \\
! Inputs: \(A^{(e)}=S^{(e)}+M^{(e)}\), \(A^{(n)}\), \(D_{h}^{(n\rightarrow e)}\) \\
! Outputs: \(P^{(e)}\), \(D_{H}^{(n\rightarrow e)}\) \\
Form \(P^{(n)}\) from \(A^{(n)}\) using AMG algorithm of choice for Poisson problem \eqref{eq: nodal problem}\\
\(\overline{P}^{(n)}\leftarrow \textsc{Gen\_Piecewise}(P^{(n)},\text{nPasses})\) (Algorithm~\ref{algo: constant interpolation}) \\
Build \(D_{H}^{(n\rightarrow e)}\) using \(D^{(n\rightarrow e)}_{h}\), \(\overline{P}^{(n)}\) \\
$[~{\cal N}, Q, R, D^{(n\rightarrow e)}_{H}~ ]$ = \textsc{Emin\_Setup}($~D^{(n\rightarrow e)}_h,D^{(n\rightarrow e)}_{H},P^{(n)}~$) (Algorithm~\ref{algo: setup emin})  \\
Solve energy minimization problem \eqref{eq:min2_hcurl} for \(P^{(e)}\) using \({\cal N}\), \(P^{(n)}\), \(D^{(n\rightarrow e)}_{h}\), \(D^{(n\rightarrow e)}_{H}\)
\end{tabbing}
\caption{\label{algo: construction of edge prolongator} Pseudo-code to construct $P^{(e)}$ and \(D^{(n\rightarrow e)}_{H}\)}
\end{algorithm}

\section{Ideal Behavior in Two Dimensions}\label{sec:ideal}
In certain limited scenarios, the SpHcurlAMG algorithm produces an ideal interpolation operator when the supplied nodal prolongator $P^{(n)}$
is based on linear interpolation and $\sigma$ is constant. In particular, the columns of the resulting edge interpolation operator  $P^{(e)}$
correspond to basis functions associated with linear edge elements in some 2D cases. This occurs for triangular-shaped
elements where a fine mesh is obtained by refining the coarse mesh elements. This includes cases where the refinement factor is
greater than two and when the refinement is not necessarily uniform within the elements. Something similar occurs for quadratic elements
on tensor product
meshes. We briefly explain the intuition behind this for meshes based on triangular elements (similar arguments hold for
tensor product meshes). We defer a detailed proof to a future monograph. %on the subject.

Consider a triangle element coarse mesh tessellation ${\cal G}_H$ and a refinement process that produces a fine mesh ${\cal G}_h$.
Discretization of $\nabla \times \nabla \times $ on the fine mesh gives rise to the matrix $S^{(h)}_h$.
% where we now omit the superscript $(e)$ and the subscript $h$ to avoid cluttering in the presentation below.
We restrict ourselves
to the case where natural boundary conditions are prescribed on all domain boundaries to again simplify the discussion.
Let the function $\phi$ be a first order coarse edge basis function that has nonzero support
within two triangles or just one depending or not whether its peak lies along an interior coarse edge or a  domain boundary edge. %is if it lies on the domain boundary.   %  This basis function support will lie within
%these adjacent elements.
As the basis function is linear within each coarse element interior,
it follows that $\nabla \times \nabla \times \phi = 0$ within the interior of each
element and that components of $w_h$ associated with interior fine degrees-of-freedom are also zero when
%$w_h = S^{(e)}_h \widetilde{\phi}$. Here, we introduce a tilde accent to represent a discrete version of $\phi$ in lieu of an $h$ subscript to avoid multiple subscripts in what follows.
$w_h = S^{(e)}_h \phi_h$. Here, $\phi_h$ represents a discrete version of $\phi$. %in lieu of an $h$ subscript to avoid multiple subscripts in what follows.

Let us now define an {\it ideal} edge interpolation operator whose columns are given by the
different $\phi_h$ basis functions.  It is well-known
that compatible discretizations are structure preserving in that they yield discrete operators that
preserve the continuous null-space property and %that
an analog of our commuting relationship involving
transfers between the continuous space and subspaces defining the discrete representation. It is also known (and can be inferred
from this structure preserving feature) that our commuting relationship is also satisfied when edge interpolation operators are based on
standard finite element edge basis functions and nodal interpolation operators are based on standard nodal elements.
That is, the ideal edge interpolation operator is a feasible solution to our energy minimization process when
the nodal prolongator is based on standard nodal elements.
Let us further suppose that the {\it ideal} edge interpolation operator is provided
as an initial guess to the energy minimization process whose version of the saddle point system given in \eqref{eq: KKT}
is
\begin{equation} \label{eq: KKT2}
\left [
\begin{matrix}
\widehat{S}^{(e)}_h  & X^T \\
X  & 0
\end{matrix}
\right ]~\left [
\begin{matrix}
\widebreve{\delta P^{(e)}} \\
\lambda
\end{matrix}
\right ]
=
\left [
\begin{matrix}
-\widebreve{S^{(e)}_h P_0^{(e)}} \\
0
\end{matrix}
\right ]
\end{equation}
where now $\hat{S}^{(e)}_h$ is a block matrix with $S^{(e)}_h$ repeated multiple times along the block diagonal
and
$P^{(e)}_0$ is the ideal edge interpolation operator that serves as a feasible initial guess.
The matrix $X$ again represents the two types of constraints
and $\delta P^{(e)}$ is a correction to $P^{(e)}_0$ giving the desired optimal solution $P^{(e)} \leftarrow P^{(e)}_0 + \delta P^{(e)}$.
Thus, our goal is to show that $\delta P^{(e)} = 0$ satisfies the above linear system.

Notice that the gradient term $-S^{(e)}_h P^{(e)}_0$ must be zero for all rows associated with fine grid
edges that reside within the interior of the coarse triangles. This follows from the
above $S^{(e)}_h \phi_h$ discussion. Thus, the only nonzero right hand side components of
\eqref{eq: KKT2} are associated with fine edges that lie along a coarse edge of ${\cal G}_H$. We recall, however,
that the prolongator pattern for a row associated with one of these fine edges has only one nonzero.
This implies that certain rows and columns of $X$ have only one nonzero if that row of $X$ enforces a
sparsity constraint or a commuting constraints for a fine edge that lies on a coarse edge.
We can then infer
that there must trivially exist a $\lambda$ such that
the only nonzero elements of $ X^T \lambda$  are equal to the only nonzero entries in the vector associated with $S^{(e)}_h P^{(e)}_0$
as these all correspond to cases where $X$ has only one nonzero in both its row and column.
%associated with fine edges that lie on coarse edges are equal
It then follows that $\delta P^{(e)} = 0$ is a solution to~\eqref{eq: KKT}.
It is important to recognize that this solution of the saddle-point system is not necessarily unique
due to the singularity of $S^{(e)}_h$. That is, the ideal prolongator is not necessarily the only solution that might
be determined by an energy minimization process. We have observed that this ideal solution is
in fact found on uniformly refined (with refinement factors that can be greater than two) triangular
meshes and for refined tensor product meshes for quadrilaterals. We observe that the energy minimization process
applied to prolongators associated with irregular refinements are {\it close} but not identically equal to ideal
interpolation operators. This ideal property does not hold in three dimensions as fine edges that lie on
faces between coarse elements do not always lie on a coarse edge. This means that while  $S^{(e)}_h P^{(e)}_0$ is zero within element
interiors, the nonzeros components do not generally correspond to rows where the sparsity prolongator
has only one nonzero and so the arguments concerning $ X^T \lambda$ do not hold. Simply put, the
energy minimization process will attempt to smooth coarse basis functions at the peaks (where
the ideal functions have weak discontinuities or kinks) along fine edges that lie on faces
where two coarse elements meet.

Finally, we note that in practice we have not observed any numerical issues associated with the singularity
of $S^{(e)}_h$. % and the above saddle-point system.
Normally,  the basis functions for the energy minimization initial guess are smooth
%does not have large components within the null space of $S^{(e)}_h$
and a few projected-Jacobi or projected conjugate gradient sweeps tends not to introduce oscillatory components
as
%these null space components are not significantly amplified by
these sweeps involve $S^{(e)}_h P^{(e)}$ computations that obviously
include no components in the null space of $S^{(e)}_h$. We have considered the addition of perturbations (gauge terms) to $S^{(e)}_h$
that guarantee the non-singularity of the $(1,1)$ block of the saddle point system, but have seen no computational benefit
from including these terms.

%%%%%%%%%%%%%%%%%%%%%%%%%%%%%%%%%%%%%%%%%%%%%%%%%%%%%%%%%%%%%%%%%%%%%
%%%%%%%%%%%%%%%%%%%%%%%%%%%%%%%%%%%%%%%%%%%%%%%%%%%%%%%%%%%%%%%%%%%%%
%%%%%%%%%%%%%%%%%%%%%%%%%%%%%%%%%%%%%%%%%%%%%%%%%%%%%%%%%%%%%%%%%%%%%
\section{Computational Examples}\label{sec:examples}
\subsection{Model Problems}
We first consider the solution of model problem~\eqref{eq: eddy current} with constant $\sigma $ and natural boundary conditions
on a square or cuboid mesh in two and three dimensions respectively.
First order
edge elements are used to discretize the PDE in each case. For the two dimensional problems we consider both uniform quadrilateral meshes and triangular meshes. The triangular meshes are obtained by
splitting each quadrilateral into two triangles. For the three dimensional meshes we consider both tetrahedral and hexahedral meshes. Tetrahedrons are created by splitting each hexahedron into six tetrahedrons.
Table~\ref{tab: model2d} gives the conjugate gradient iterations for the two dimensional cases
\begin{table}[h!]
\centering
\caption{CG iterations and AMG operator complexity (AMG \textsf{o.c.}) for 2D constant $\sigma$ model problems.}
\label{tab: model2d}
\begin{tabular}{c|c|ccccc|c}
  \multicolumn{7}{c}{triangular mesh} \\[5pt]
\toprule
nodal mesh     & \#edges  &  \multicolumn{5}{c|}{$\sigma$ } & AMG \textsf{o.c.} \\[3pt]
               &           &     $10^2$        &    $10^{1}$   &    $10^{0}$   &    $10^{-1}$  &     $10^{-2}$    \\
$ 28\times  28$&$    2,241$&       3        &       8        &    9           &     9          &    ~$8 $~       & 1.17\\
$ 82\times  82$&$   19,845$&       7        &       9        &    9           &    10          &    ~$8 $~       & 1.20\\
$244\times 244$&$  177,633$&       8        &       9        &   10           &     9          &    ~$7 $~       & 1.19\\
$730\times 730$&$1,595,781$&       9        &       9        &    9           &     7          &    ~$5 $~       & 1.19\\
% tolerance changes (top to bottom) mu=1e 0:   -    ,     -   ,     -   ,  1.68e-8
% tolerance changes (top to bottom) mu=1e-1:   -    ,  1.53e-8,  1.65e-7,  1.69e-6
% tolerance changes (top to bottom) mu=1e-2: 1.19e-7,  1.32e-6,  1.65e-5,  1.69e-4
    \bottomrule
    \end{tabular}
\vskip .2in
\begin{tabular}{c|c|ccccc|c}
  \multicolumn{7}{c}{quadrilateral mesh} \\[5pt]
\toprule
nodal mesh     & \#edges  &  \multicolumn{5}{c|}{$\sigma$ }  & AMG \textsf{o.c.} \\[3pt]
               &           &     $10^2$        &    $10^{1}$   &    $10^{0}$   &    $10^{-1}$  &     $10^{-2}$    \\
$ 28\times  28$&$    1,512$&      2         &     5          &    6           &      6         &    ~$6$ & 1.11 \\
$ 82\times  82$&$   13,284$&      4         &     6          &    6           &      6         &    ~$5$ & 1.13 \\
$244\times 244$&$  118,584$&      4         &     6          &    6           &      6         &    ~$5$ & 1.13 \\
$730\times 730$&$1,064,340$&      6         &     6          &    6           &      5         &    ~$4$& 1.13 \\
% tolerance changes (top to bottom) mu=1e-1:   -    ,     0   ,  2.76e-8,  2.49e-7
% tolerance changes (top to bottom) mu=1e-2: 3.34e-8,  4.64e-7,  2.76e-6,  2.49e-5
    \bottomrule
    \end{tabular}
  \end{table}
using one AMG V-cycle sweep as a preconditioner.  The table also shows
the AMG operator complexity \textsf{o.c.}, 
which is the ratio between the sum obtained by adding the number of discretization matrix nonzeros
on all levels and dividing this by the number of nonzeros on just the finest level. 
The nodal AMG hierarchy is constructed by applying a standard smoothed
aggregation multigrid method to the nodal problem~\eqref{eq: nodal problem} where the drop tolerance is zero. Only one Jacobi-style iteration of the energy minimization procedure is used to generate
edge prolongators using $\omega = .5$ in \eqref{emin iteration}. The initial prolongator {\it guess} is obtained by setting nonzeros in the sparsity pattern to one and then solving
local least squares problems to correct the prolongator and satisfy the commutator equation. For these model cases, the
Reitzinger and Sch\"oberl algorithm~\eqref{eq: revised RS DH} is sufficient for generating the $D^{(n\rightarrow e)}_H $ matrices without requiring additional edges associated
with make\_irreducible() mentioned in Algorithm~\ref{algo: setup emin}. As noted, Algorithm~\ref{algo: constant interpolation} is used to generate the piecewise constant interpolation
operators needed for the $D^{(n\rightarrow e)}_H $ construction.

In Table~\ref{tab: model2d}, a V-cycle is employed with 
one pre- and one post-symmetric Hiptmair iteration 
for AMG relaxation~\cite{Hiptmair98}.
On the finest level, one symmetric Hiptmair sweep corresponds to applying one symmetric Gauss-Seidel iteration to $A_e u_h = f_h$. This is then followed by performing
one symmetric Gauss-Seidel iteration applied to the system $[D^{(n\rightarrow e)}_h]^T  A_e D^{(n\rightarrow e)}_h c_h = [D^{(n\rightarrow e)}_h]^T r_h$ where $c_h$ is initialized to zero and
$r_h$ is the residual
after the initial symmetric Gauss-Seidel iteration. The vector $D^{(n\rightarrow e)}_h c_h$ is then used to correct the solution. Finally, one additional
symmetric Gauss-Seidel iteration is applied to the $A_e$ system. Hiptmair relaxation is applied on all levels except the coarsest where a direct solver
is employed. The number of AMG levels ranges from two for the smallest problems to five for the finest problems.
The initial guess to CG is zero and the right hand side is chosen as a random vector.
The default CG convergence tolerance is $ \tau = 10^{-8} ||b||$ where $||b||$ is the 2-norm of the right hand side. Due to the potentially large condition numbers (for small $\sigma$), we find that rounding
errors can play a role in the achievable accuracy (even for a direct solver). To account for this, we apply Matlab's \textsf{condest()} function to the matrix obtained by left scaling
the finest level discretization operator using the inverse matrix diagonal to define the scaling. This condition number estimate $\kappa(D^{-1} A)$  is then used to change the convergence tolerance
so that it does not exceed $\kappa(D^{-1} A)/10^{16}$. 
This slightly changes the tolerance for $\sigma = 1.0$ on the finest triangular mesh. It changes the $\sigma = 10^{-1}$ triangle element
tolerances by approximately a factor of $15$ and $150$ the $2^{nd}$ finest and finest meshes respectively. 
The $\sigma = 10^{-2}$ triangle      element tolerances are changed by approximately a factor of $12,130,1650,$ or $17000$ from the coarsest to the finest mesh respectively.
The $\sigma = 10^{-1}$ quadrilateral element tolerances are changed by approximately a factor of $3$ and $25$ on the $2^{nd}$ finest and finest mesh respectively.
The $\sigma = 10^{-2}$ quadrilateral element tolerances are changed by approximately a factor of $3, 50, 275$, and $2500$ from the coarsest to the finest mesh respectively.
The associated largest value for  $\kappa(D^{-1} A)$ is $1.7 \times 10^{12}$, which occurs for the finest triangle mesh problem with $\sigma = 10^{-2}$.

As can be seen in Table~\ref{tab: model2d}, the iteration count is scalable and the AMG operator complexity
is only slightly greater than one. Additionally, there is only a modest convergence sensitivity to the value of $\sigma$ (with modest reductions in iterations due to modified tolerances).
Finally, we note that these AMG iteration counts
are comparable to those obtained by applying geometric multigrid to the same problem, which requires $8$ iterations or $6$ iterations for the triangular or
quadrilateral meshes respectively when $\sigma = 10^{-2}$. In the geometric multigrid case, a coarsening rate of three is used (e.g., the $244 \times 244$ problem uses a $82 \times 82$ and a
$28 \times 28$ coarse mesh) and the V-cycle particulars are the same as for the AMG method.

Table~\ref{tab: 2d least squares} provides statistics on the dimensions of the least-squares computations required to enforce the
\begin{table}[h!]
\centering
\caption{Number of least squares problems with corresponding matrix dimension during AMG setup phase for 2D model problem. }
\label{tab: 2d least squares}
\begin{tabular}{c|ccccc}
  \multicolumn{6}{c}{triangular mesh with $1595781$ edges} \\[5pt]
\toprule
               &   \multicolumn{5}{c}{least squares sub-problem dimensions } \\
               &  ~~$1 \times 2$~~ & ~~$2 \times 3$~~ & ~~$3 \times 3$~~ & ~~$5 \times 4$~~ & ~~$6 \times 4$~~ \\ \hline & & & & & \\[-6pt]
$  P^{(e)}_0 $ &  $   355752 $& $    930   $ & $   117168 $ & $ 1112754  $ & $ 9177  $ \\
$  P^{(e)}_1 $ &  $   39852  $& $    165   $ & $    13120 $ & $  121639  $ & $ 3340  $ \\
$  P^{(e)}_2 $ &  $    4536  $& $     53   $ & $     1458 $ & $   12885  $ & $ 1072  $ \\
$  P^{(e)}_3 $ &  $     540  $& $     32   $ & $      147 $ & $    1551  $ & $   22  $ \\
    \bottomrule
    \end{tabular}
\vskip .2in
\begin{tabular}{c|ccccc}
  \multicolumn{6}{c}{quadrilateral mesh with $1064340$ edges} \\[5pt]
\toprule
               &   \multicolumn{5}{c}{least squares sub-problem dimensions } \\
               &  ~~$1 \times 2$~~ & ~~$2 \times 3$~~ & ~~$3 \times 3$~~ & ~~$4 \times 4$~~ & ~~$5 \times 4$~~ \\ \\ \hline & & & & & \\[-6pt]
$  P^{(e)}_0 $ &  $   355752 $& $      0   $ & $      0   $ & $ 697212   $ & $ 11376  $ \\
$  P^{(e)}_1 $ &  $   39852  $& $      0   $ & $     320  $ & $ 74892    $ & $ 4468  $ \\
$  P^{(e)}_2 $ &  $    4536  $& $      0   $ & $     106  $ & $ 7500     $ & $ 1462  $ \\
$  P^{(e)}_3 $ &  $     540  $& $     30   $ & $       4  $ & $  995     $ & $   47  $ \\
    \bottomrule
    \end{tabular}
  \end{table}
commutator equation during the AMG setup phase for the highest resolution triangle-mesh and quadrilateral-mesh problems. Here, $P^{(e)}_\ell$
denotes the prolongator setup on level $\ell$ in the multigrid hierarchy (with $\ell = 0$ being the finest). We can see that the largest least
squares problems are only $ 6 \times 4$ on the triangular mesh while the vast majority of the sub-problems are $5 \times 4$ or smaller.
As there are no Dirichlet boundary conditions, the rank is one lower than the column dimension. Thus, using the normal equations for a
$5 \times 4$ system would only require the inversion of a $3 \times 3$ matrix.  It should be noted that the somewhat limited number of least squares problem sizes
and their small size is due in part to the regularity of our aggregation algorithms when applied to a structured mesh.  Interestingly,
many of the matrices with the same dimension must actually be identical up to a permutation or sign changes as they all represent
a sub-graph connecting a small set of nodes.
More irregular
meshes or irregular coarsening would lead to a wider range of least squares sub-systems.  As mentioned earlier, the column dimensions
of the linear algebra sub-systems to solve 3D elasticity problems would normally be $ 6 $ (and would be $3$ in two dimensions).
Thus, loosely our linear sub-systems are
of comparable size to those required for elasticity (perhaps a bit larger than required for 2D elasticity but a bit smaller
than required for 3D elasticity). While our current Matlab implementation do not provide credible timings, the comparable linear
systems sizes give us some hope/expectation that the associated setup time of a proper implementation will not be long
compared to the solve time as was demonstrated in~\cite{JaFrScOl2023} for elasticity.

Table~\ref{tab: model3d} gives the results for the three dimensional cases.
\begin{table}[h!]
\centering
\caption{CG iterations and AMG operator complexity (AMG \textsf{o.c.}) for 3D constant $\sigma$ model problems.}
\label{tab: model3d}
\begin{tabular}{c|c|ccccc|c}
  \multicolumn{7}{c}{tetrahedral mesh} \\[5pt]
\toprule
nodal mesh     & \#edges  &  \multicolumn{5}{c|}{$\sigma$ } & AMG \textsf{o.c.} \\[3pt]
               &            &     $10^2$             &    $10^{1}$   &    $10^{0}$   &    $10^{-1}$  &     $10^{-2}$    \\
$ 10\times  10\times 10$&$    5,859$&  4             &   9            &   11           &    12          &    ~$12$  & 1.13 \\
$ 28\times  28\times 28$&$  144,423$&  5             &  12            &   12           &    13          &    ~$13$  & 1.11 \\
$ 82\times  82\times 82$&$3,779,379$&  9             &  12            &   13           &    13          &    ~$11$  & 1.10 \\
% tolerance changes (top to bottom) mu=1e-1:   -    ,     -   ,  2.58e-8
% tolerance changes (top to bottom) mu=1e-2: 3.18e-8,  2.87e-7,  2.58e-6
    \bottomrule
    \end{tabular}
\vskip .2in
\begin{tabular}{c|c|ccccc|c}
  \multicolumn{7}{c}{hexahedral mesh} \\[5pt]
\toprule
nodal mesh     & \#edges  &  \multicolumn{5}{c|}{$\sigma$ } & AMG \textsf{o.c.} \\[3pt]
               &                   &     $10^2$        &    $10^{1}$   &    $10^{0}$   &    $10^{-1}$  &     $10^{-2}$    \\
$ 10\times  10\times 10$&$    2,700$&    3           &   4            &    6           &     6          &       6           & 1.09 \\
$ 28\times  28\times 28$&$   63,504$&    3           &   5            &    6           &     6          &       5           & 1.06 \\
$ 82\times  82\times 82$&$1,633,932$&    3           &   6            &    6           &     6          &       5           & 1.05 \\
% tolerance changes (top to bottom) mu=1e-1:   -    ,     -   ,  1.85e-8
% tolerance changes (top to bottom) mu=1e-2: 2.41e-8,  2.04e-7,  1.85e-6
    \bottomrule
    \end{tabular}
  \end{table}
All algorithm details are identical to the two dimensional cases with similar modified tolerances. 
The $\sigma = 10^{-1}$ tetrahedral element tolerance is changed by approximately a factor of $2.5$ on the finest mesh.
The $\sigma = 10^{-2}$ tetrahedral   element tolerances are changed by approximately a factor of $3,30,$ and $260$ from the coarsest to the finest mesh respectively.
The $\sigma = 10^{-1}$ hexahedral element tolerance is changed by approximately a factor of $2$ on the finest mesh.
The $\sigma = 10^{-2}$ hexahedral    element tolerances are changed by approximately a factor of $   2.5, 20,$ and $185$ from the coarsest to the finest mesh respectively.
We note that using two energy minimization iterations to compute the edge prolongators does reduce the iteration
count by one for a few of the cases considered in the table and never increases the iteration count for any of these cases.
For this problem, geometric multigrid requires $11, 12,$ and $13$ iterations as one refines the tetrahedral mesh
respectively for the smallest $\sigma$ problem, which is comparable to that obtained by AMG. 
Geometric multigrid
experiments where not run for the hexahedral mesh as we do not currently have this capability.

Table~\ref{tab: 3d least squares} provides the least-squares statistics for the 3D case.
\begin{table}[h!]
\centering
\caption{Number of least squares problems with corresponding matrix dimension during AMG setup phase for 3D model problem. }
\label{tab: 3d least squares}
\begin{tabular}{c|cccccc}
  \multicolumn{7}{c}{tetrahedral mesh with $3779379$ edges} \\[5pt]
\toprule
\z\z\z            \z\z\z    &   \multicolumn{5}{c}{least squares sub-problem dimensions } \\
               &  ~$1 \times 2  $& $[2\z:\z 3] \times 3$  & $[5\z:\z 6] \times 4$ & $[7\z:\z 10] \times 5$ & $[10\z:\z 15] \times 6$ &  $[19\z:\z 24] \times 8$\\ \\[-9pt] \hline & & & & & \\[-6pt]
$  P^{(e)}_0 $ &   $   190512   $& $    122472          $ & $   1163484         $ & $ 118098              $ & $354294              $ & $1830519$  \\
$  P^{(e)}_1 $ &   $    8100    $& $    4860            $ & $    46170          $ & $  4374              $ & $ 13122              $ & $ 67797 $  \\
$  P^{(e)}_2 $ &   $     432    $& $     240            $ & $    2106           $ & $   160              $ & $  552               $ & $ 2660  $  \\
    \bottomrule
    \end{tabular}
\vskip .2in
\begin{tabular}{c|cccccc}
  \multicolumn{6}{c}{hexahedral mesh with $1633932$ edges} \\[5pt]
\toprule
               &   \multicolumn{6}{c}{least squares sub-problem dimensions } \\
               &  ~$1 \times 2$~ & ~$[2\z:\z  3] \times 3$~ & ~$[4\z:\z  6] \times 4$~ & ~$[6\z:\z 10] \times 5$~ & ~$[9\z:\z 13] \times 6$~ & ~$[12\z:\z19] \times 8$ \\ \\[-9pt] \hline & & & & & \\[-6pt]
$  P^{(e)}_0 $ &  $   190512 $     &                            &   $   734832  $            &                            &                            &   $ 708588    $ \\
$  P^{(e)}_1 $ &  $    8100  $     &   $   1110     $           &   $    32126  $            &   $ 1302               $   &   $ 2730               $   &   $  32469    $ \\
$  P^{(e)}_2 $ &  $     432  $     &   $    150     $           &   $     1759  $            &   $  114               $   &   $  276               $   &   $   1549    $ \\
    \bottomrule
    \end{tabular}
  \end{table}
The notation $[x\z:\z y] \times z$ indicates that the row dimension of these least-squares problems have been grouped together with the smallest
and largest row dimensions being $x$ and $y$ respectively. As these sub-problems are again singular, a normal equations approach
word require at most the solutions to $7 \times 7$ linear systems, which is still comparable to that required for 3D elasticity.
Additionally, we note that in some sense SpHcurlAMG sub-systems are {\it easier} to solve than those required for 3D elasticity.
For elasticity, the rectangular matrices are fully dense while in our context they are sparse (two $\pm 1$ nonzeros per row).
This implies that less storage is needed and matrix-vector operations require less flops, though
leveraging sparsity would probably only be beneficial for fairly large sub-systems. More importantly, the rank is known in advance,
which is not the case in elasticity where one must address the conditioning of the least-squares problems.

A full assessment of the setup times will be performed once we have completed a proper algorithm implementation. In lieu of this
assessment, we make some additional remarks regarding these least-squares systems. Ultimately, the nodal sparsity of the nodal
prolongator determines the column dimensions of the sub-systems. As is well-known, smoothed aggregation (used for these
experiments) typically generates relatively sparse interpolation operators leading to low operator complexity. If other AMG algorithms
are employed to generate the $P^{(n)}$'s, it might be important in our context to consider variants that
lead to relatively sparse interpolation operators such as those in ~\cite{de_sterck_distance-two_2008}.
We also note that there are potential short-cuts that one could consider to reduce setup time. This includes
leveraging at least part of the setup from a previous solve when a sequence of linear systems is solved (e.g., for implicit time stepping).
For example, the interpolation operators can often be retained across a few linear solves
(still re-projecting the discretization matrix for each solve). Another possibility might be to retain the
nodal AMG strength-of-connection matrix so that the
least-square factors can be retained as the sparsity patterns of $P^{(n)}$ and $P^{(e)}$ do not change.
Finally, we mention that other techniques (related to randomized linear algebra) could be considered
to further reduce costs. In our context, the sub-systems are under-determined and so the QR or normal
equations is employed to obtain a least-norm solution. One could instead approximate the least-norm solution by combining solutions
to a few square systems (which can each be chosen to be lower triangular).

\subsection{ALEGRA Problems}
\begin{figure}[htb!]\label{fig: badapple}
%\hskip -0.05in {\includegraphics[trim=240 0 240 0,clip,height=2.9in,width=1.71in]{Figs/apple_6.png}} \vskip -2.57in \hskip  1.75in
 \hskip -0.05in {\includegraphics[trim=240 0 240 0,clip,height=2.9in,width=1.71in]{     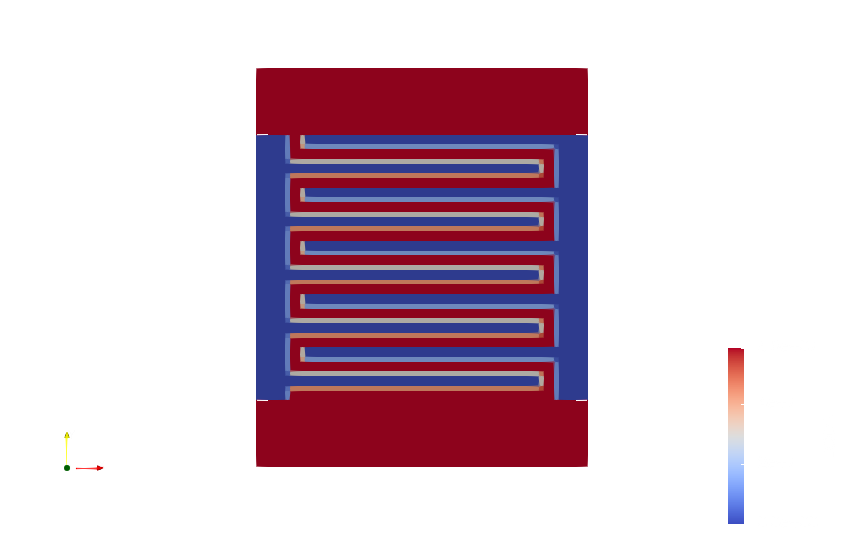}} \vskip -2.57in \hskip  1.75in
%{\includegraphics[trim=45 210 45 210,clip,height=2.2in,width=3.1634in]{Figs/wideapple6.jpg}}
 {\includegraphics[trim=45 210 45 210,clip,height=2.2in,width=3.1634in]{     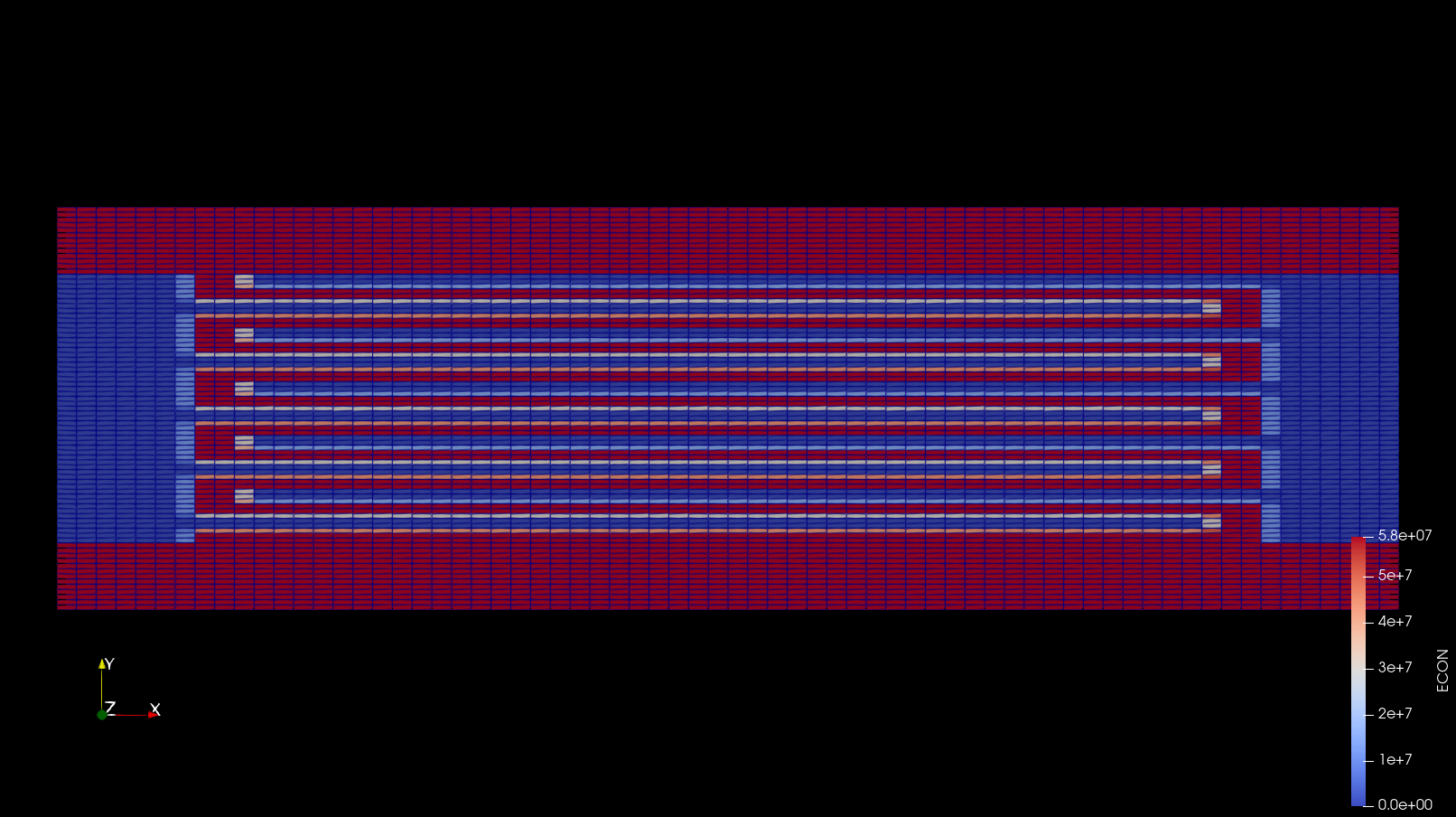}}
\caption{Snake-like conductor problem. Left image shows material variation for the {\it standard} problem version. Right image shows the hexahedral mesh for the {\it wide} wider version.}
\end{figure}

Figure~\ref{fig: badapple} illustrates an ALEGRA~\cite{AlegraALE08,ROBINSON2024117164,NIEDERHAUS2023104693} simulation test problem that employs first order edge elements on a regular hexahedral mesh. The mesh
for the finest resolution ``standard'' problem is $ 69 \times 82 \times 3$ while the ``wide'' version employs a $ 82 \times 69 \times 3$ mesh.
The long
snake-like red conductor poses some AMG challenges in that the coarsening must generally preserve the snake-like structure
throughout the AMG hierarchy.
Table~\ref{tab: badapple} shows results using a 2-level AMG method on a snake-like conductor matrix which has been obtained by first exporting it from the ALEGRA code (i.e., saving it to disk).
The nodal prolongator was produced using the Trilinos/MueLu code
by applying smoothed aggregation to the associated nodal problem. This method is fairly standard except that a drop tolerance of $.0025$ is used in conjunction with a
standard symmetric smoothed aggregation dropping criteria. The drop tolerance attempts to form aggregates that respect material boundaries within the smoothed
aggregation algorithm.
While these problems are relatively small, an estimate of the matrix condition number
is approximately $4.5 \times 10^{12}$ due to the more than 7 orders of magnitude in the  material variation. %This estimate was obtained
\begin{table}[h!]
\centering

\caption{AMG applied to snake conductor problem. }
\label{tab: badapple}
\begin{tabular}{c|c| c|c|c}
problem &   \# edges                      & precond & \#CG \textsf{its.}&  AMG \textsf{o.c.}\\
\toprule
\multirow{6}{15mm}{\\standard} &\multirow{2}{12mm}{19857}                                           & Relaxation only  & $  33   $& $     1.0   $   \\[3pt]
                         &                                                                    & 2-level AMG      & $   6   $& $     1.02  $   \\[2pt] \cmidrule{2-5} & & & & \\[-6pt]
                         &\multirow{2}{12mm}{31715}                                           & Relaxation only  & $  64^* $& $     1.0   $   \\[3pt]
                         &                                                                    &  2-level AMG     & $   7   $& $     1.02  $   \\[2pt] \cmidrule{2-5} & & & & \\[-6pt]
                         &\multirow{2}{12mm}{43573}                                           & Relaxation only  & $  78^*  $& $     1.0   $   \\[3pt]
                         &                                                                    &  2-level AMG     & $   7   $& $     1.02  $   \\[2pt] \hline& & & & \\[-6pt]
\multirow{2}{15mm}{wide} &\multirow{2}{12mm}{43573}                                           & Relaxation only  & $  78^* $& $     1.0   $   \\[3pt]
                         &                                                                    &  2-level AMG     & $   8   $& $     1.02  $   \\[2pt]
    \bottomrule
    \end{tabular}
  \end{table}
The column dimension of the largest least squares sub-problem is $6$, while the vast majority of these least square problems have column dimension $2$ or $4$.
                        For these results, CG is used as the outer Krylov method with a tolerance of approximately $9\times 10^{-6} $ while other method details are the same as in the model problem case.
The $9\times 10^{-6} $ tolerance is computed as $\kappa(D^{-1}A)/(5 \times 10^{17})$, which was experimentally determined as a reasonable threshold for this matrix. %We find that a direct
An asterisk for the relaxation only solves
indicates that conjugate gradient failed to achieve the tolerance, but came fairly close (within a factor of 10). Specifically,
the CG method stagnated with two consecutive iterations with the same solution. It it possible to achieve the tolerance by effectively restarting CG and taking a
few additional iterations (which are not included in the table). As noted, this is due to the extreme ill-conditioning in this problem.
For these experiments the Jacobi energy minimization procedure uses a damping parameter of $.3$ and only one energy minimization iteration is taken.
As the results illustrate, the AMG iterations required for convergence do not increase significantly as the problem is refined and that at most eight
iterations is needed despite the large condition number.

\begin{figure}[htb!]\label{fig: helix}
\centering
%{\includegraphics[trim=185 100 300 0,clip,scale=.3]{Figs/for_ray.png}}
 {\includegraphics[trim=185 100 300 0,clip,scale=.3]{     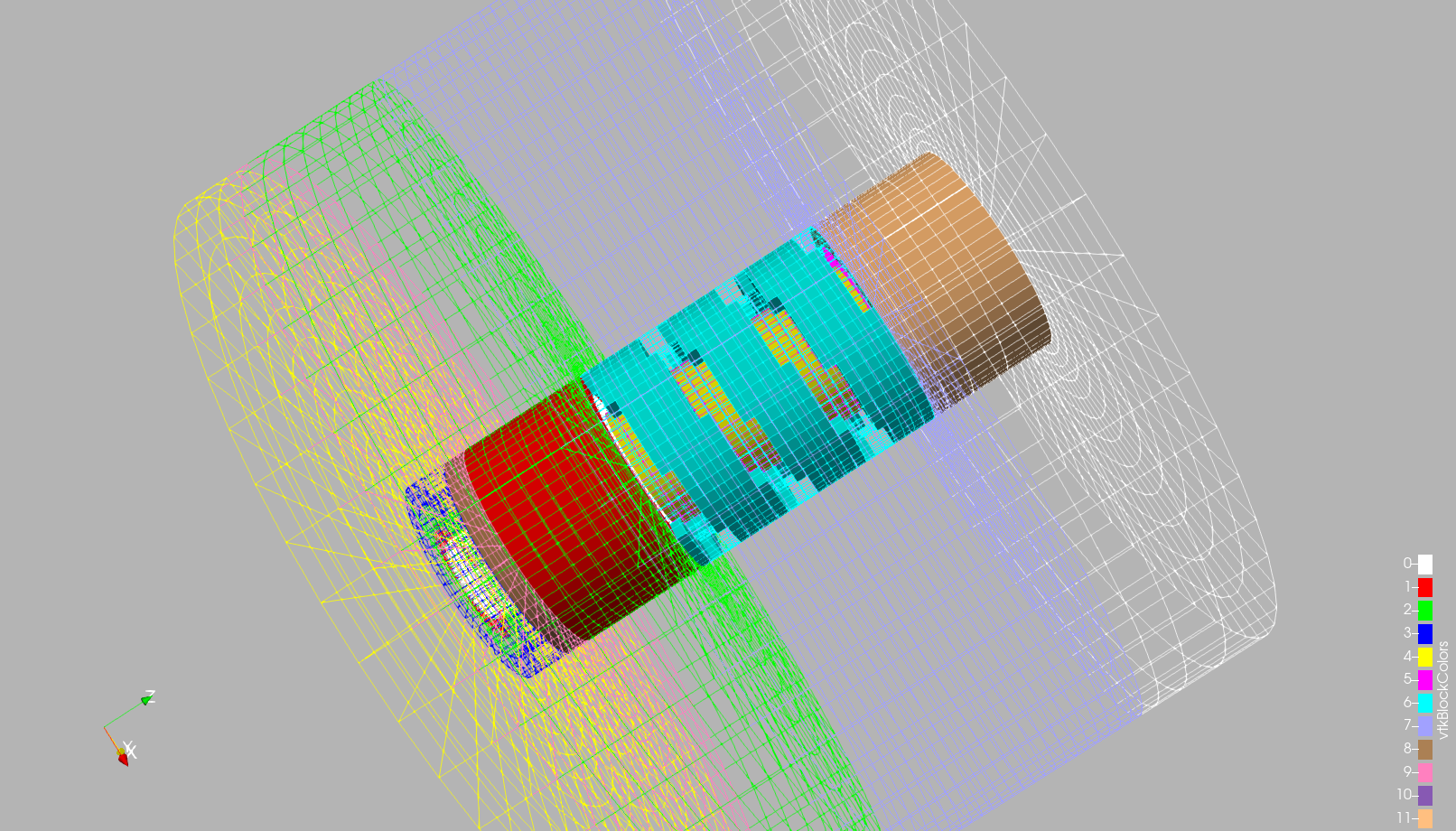}}
\caption{Coax cable mesh and material variations.}
\end{figure}

As a third test case, we consider a coax cable problem obtained by exporting a matrix from the ALEGRA simulation code.  %~\cite{AlegraALE08}.
The overall simulation involves transient magnetics within a coax cable which attempts to force helical current flow (see Figure~\ref{fig: helix}).
The large material discontinuities (over seven orders of magnitude) along with the long helical path that the current travels give rise to severely ill-conditioned linear systems. Table~\ref{tab: helix} illustrates AMG convergence rates associated with using 3 AMG levels or just 2 sweeps of the Hiptmair relaxation method (i.e., without employing multigrid) for the first linear system solve in the simulation.
Two Hiptmair sweeps is equivalent to one pre- and one post- relaxation sweep on the finest level of a 3-level method.
\begin{table}[h!]
\centering

\caption{AMG applied to coax problem. }
\label{tab: helix}
\begin{tabular}{cc| c|c|c}
\# edge dofs & $\kappa(D^{-1}A)$  & precond & \#CG \textsf{its.}&  AMG \textsf{o.c.}\\
\toprule
\multirow{2}{12mm}{164259} & \multirow{2}{20mm}{$5.9 \times 10^{11}$}& Relaxation only  &    84    &       1.0       \\[3pt]
                           &                                         & 3-level AMG      &    28    &       1.18      \\[2pt] \hline & & & & \\[-6pt]
\multirow{2}{12mm}{464715}  & \multirow{2}{20mm}{$5.9 \times 10^{11}$}& Relaxation only &   202    &       1.0       \\[3pt]
                           &                                         & 3-level AMG      &    31    &       1.13      \\[2pt] \hline
    \bottomrule
    \end{tabular}
  \end{table}
The table also provides a condition number estimate for the diagonally scaled matrix. % denoted by $\kappa(D^{-1}A)$.
The nodal prolongator was again produced using the Trilinos/MueLu code with a drop tolerance of $.04$.
For these results, CG is again used as the outer Krylov method with a tolerance of $1.77 \times 10^{-6}$ while other method details are the same as in the model problem case.
As with the previous example, the tolerance is simply $\kappa(D^{-1}A)/(5 \times 10^{17})$, which was experimentally determined as a reasonable threshold for this matrix.
The initial guess is zero and the right hand side is chosen as a random vector.
As the results show, AMG reduces the iteration count by nearly a factor of $3\times$ for the smaller problem compared to the method that does not employ AMG while the reduction
is greater than $6\times$ for the larger problem. We also see that while the iteration count is not completely constant, the growth in iterations is relatively modest.

\begin{figure}[htb!]\label{fig: helix ls}
\centering
%{\includegraphics[trim=5 5 0 5,clip,height=2in,width=2.5in]{Figs/large2HelixPe0.jpg}}~~
 {\includegraphics[trim=5 5 0 5,clip,height=2in,width=2.5in]{     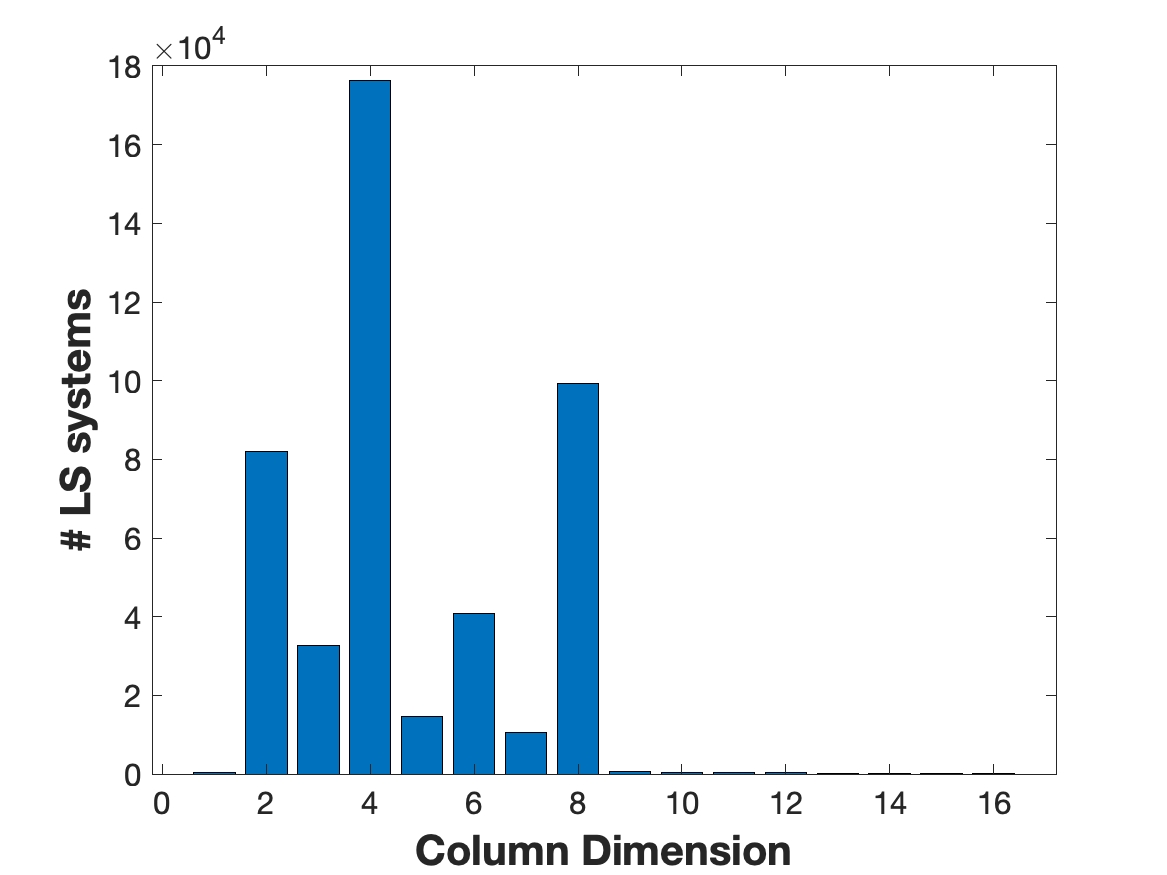}}~~
%{\includegraphics[trim=5 5 0 5,clip,height=2in,width=2.5in]{Figs/large2HelixPe1.jpg}}
 {\includegraphics[trim=5 5 0 5,clip,height=2in,width=2.5in]{     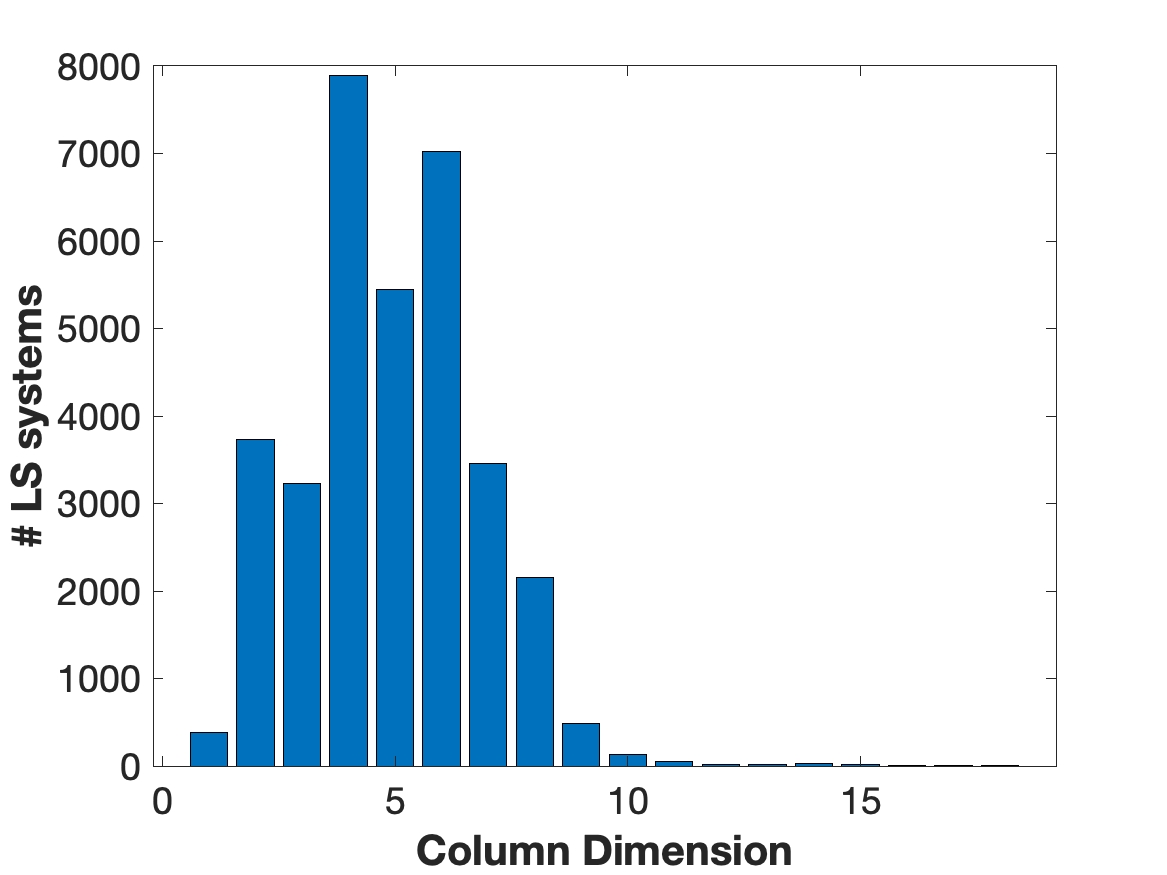}}
\caption{Number of least squares sub-systems with a particular column dimension for large coax cable case. Finest (coarsest) prolongator shown on left (right).}
\end{figure}

Figure~\ref{fig: helix ls} provides some statistics concerning the column dimension on the sub-matrices associated with
the level 0 (left) and level 1 (right) prolongators.
In both cases, the vast majority of least squares systems have column dimension of eight or less, which again is
similar to that for the 3D model problems.

%%%%%%%%%%%%%%%%%%%%%%%%%%%%%%%%%%%%%%%%%%%%%%%%%%%%%%%%%%%%%%%%%%%%%
%%%%%%%%%%%%%%%%%%%%%%%%%%%%%%%%%%%%%%%%%%%%%%%%%%%%%%%%%%%%%%%%%%%%%
%%%%%%%%%%%%%%%%%%%%%%%%%%%%%%%%%%%%%%%%%%%%%%%%%%%%%%%%%%%%%%%%%%%%%
\section{Conclusions and Future Work}\label{sec:conclusions}

We introduced a new AMG algorithm for solving the eddy current equations. 
The key feature of the new algorithm is that it is able to satisfy
a discrete commuting relationship on all multigrid hierarchy levels 
while also approximately minimizing the energy of the interpolation 
basis functions. In this way, the AMG algorithm produces structure
preserving discretizations on all hierarchy levels in the sense
that the null space properties on the finest grid are mimicked 
on the coarse grids. We have also provided some computational
evidence illustrating the convergence advantages of the 
new method. Overall, the convergence rate is similar to that
obtained with geometric multigrid on model constant coefficient
problems on regular grids. The true advantage of the approach
is that it can be applied to unstructured grids and to highly
variable coefficient problems.

%%%%%%%%%%%%%%%%%%%%%%%%%%%%%%%%%%%%%%%%%%%%%%%%%%%%%%%%%%%%%%%%%%%%%
\bibliographystyle{siamplain}
\bibliography{hcurl,YunrongZhu}
\end{document}